\documentclass{amsart} 

\usepackage{amsmath,amsthm}     
\usepackage{amssymb}            
\usepackage{euscript}           
\usepackage{epic}
\usepackage{eepic,enumerate,calc}
\usepackage[matrix,arrow,curve,frame]{xy}    

\xymatrixcolsep{1.9pc}                          
\xymatrixrowsep{1.9pc}
\newdir{ >}{{}*!/-5pt/\dir{>}}                  

\newtheorem{prop}[subsection]{Proposition}

\newtheorem{cor}[subsection]{Corollary}
\newtheorem{lemma}[subsection]{Lemma}

\theoremstyle{definition}  

\newtheorem{remark}[subsection]{Remark}

  {\end{list}}
  
\newcommand{\dfn}{\textbf} 

\newcommand{\mdfn}[1]{\dfn{\mathversion{bold}#1}} 

\newcommand{\Smash}             {\wedge}

\newcommand{\tens}              {\otimes}               
\newcommand{\iso}               {\cong}  

\newcommand{\cat}{\EuScript}    
\newcommand{\cA}{{\cat A}}      

\newcommand{\cF}{{\cat F}}

\newcommand{\field}[1]  {\mathbb #1} 
\newcommand{\A}         {\field A}
\newcommand{\R}         {\field R}
\newcommand{\F}         {\field F}

\newcommand{\Z}         {\field Z}
\newcommand{\C}         {\field C}
\newcommand{\Q}         {\field Q}
\renewcommand{\P}         {\field P}

\DeclareMathOperator{\Spec}{Spec}
\DeclareMathOperator{\Hom}{Hom}
\DeclareMathOperator{\Ext}{Ext}

\DeclareMathOperator{\Gr}{Gr}
\DeclareMathOperator{\Br}{Br}

\DeclareMathOperator{\chara}{char}

\newcommand{\ra}{\rightarrow}                   
\newcommand{\lra}{\longrightarrow}              
\newcommand{\llra}[1]{\stackrel{#1}{\lra}}      

\newcommand{\fib}{\twoheadrightarrow}           

\newcommand{\inc}{\hookrightarrow}              

\newcommand{\blank}{-}                          
\newcommand{\und}{\underline}

\newcommand{\he}{\simeq}

\newcommand{\pt}{pt}

\newcommand{\rea}[1]{|{#1}|}             

\newcommand{\ceck}[1]{\Cech(#1)}         
\newcommand{\oceck}[1]{\Cech^{o}(#1)}    
\newcommand{\oreal}[1]{\rea{\oceck{U}}}  
\newcommand{\creal}[1]{\rea{\ceck{U}}}   

\newcommand{\Gal}{\text{Gal}}

\newcommand{\Cech}{\check{C}}
\newcommand{\CCech}{\v{C}ech\ }

\newcommand{\tH}{\tilde{H}}

\newcommand{\lb}{\langle}
\newcommand{\rb}{\rangle}
\newcommand{\inner}{( \blank,\blank )}
\newcommand{\br}[1]{\lb #1 \rb}
\newcommand{\dbr}[1]{\lb\lb #1 \rb\rb}
\newcommand{\brn}[1]{\{ #1 \}}
\newcommand{\mmod}{\text{mod\ }}
\newcommand{\Qa}{Q_{\und{a}}}
\newcommand{\Zlt}{\Z_{(2)}}
\newcommand{\CQa}{\tilde{C}}
\newcommand{\Sq}{Sq}
\newcommand{\intH}{\tH_{int}}
\newcommand{\tC}{\tilde{C}}
\newcommand{\Gi}{G{\field I}}
\newcommand{\mA}{\cA}
\newcommand{\bare}[2]{[#1_{1}|\cdots|#1_{#2}]}

\numberwithin{equation}{subsection}

\newenvironment{myequation}
  {\addtocounter{subsection}{1}\begin{eqnarray}}
  {\end{eqnarray}$\!\!$}

\newcounter{property}

\begin{document}

\title{Notes on the Milnor conjectures}

\author{Daniel Dugger}

\address{Department of Mathematics\\ University of Oregon\\ Eugene, OR
97403}

\email{ddugger@math.uoregon.edu}

\date{August 26, 2004}

\maketitle

\tableofcontents

\section*{Introduction}
\label{se:intro}

These lectures concern the two Milnor conjectures and their proofs:
from \cite{Vmilnor2}, \cite{OVV}, and \cite{M2}.  Voevodsky's proof of
the norm residue symbol conjecture---which is now eight years
old---came with an explosion of ideas.  The aim of these notes is to
make this explosion a little more accessible to topologists.  My
intention is not to give a completely rigorous treatment of this
material, but just to outline the main ideas and point the reader in
directions where he can learn more.  I've tried to make the lectures
accessible to topologists with no specialized knowledge in this area,
at least to the extent that such a person can come away with a general
sense of how homotopy theory enters into the picture.

Let me apologize for two aspects of these notes.  Foremost, they
reflect only my own limited understanding of this material.  Secondly,
I have made certain expository decisions about which parts of the
proofs to present in detail and which parts to keep in a ``black
box''---and the reader may well be disappointed in my choices.  I hope
that in spite of these shortcomings the notes are still useful.

Sections 1, 2, and 3 each depend heavily on the previous one.  Section
4 could almost be read independently of 2 and 3, except for the need of
Remark~\ref{re:Hpt}.  

\vspace{0.2in}
\section{The Milnor conjectures}

The Milnor conjectures are two purely algebraic statements in the
theory of fields, having to do with the classification of quadratic
forms.  In this section we'll review the basic theory and summarize the
conjectures.  Appendix A contains some supplementary material, where
several examples are discussed.

\subsection{Background}
Let $F$ be a field.  In some sense our goal is to completely classify
symmetric bilinear forms over $F$.  To give such a form $\inner$ on
$F^n$ is the same as giving a symmetric $n\times n$ matrix $A$, where
$a_{ij}=( e_i,e_j)$.  Two matrices $A_1$ and $A_2$ represent the
same form up to a change of basis if 
and only if $A_1=PA_2P^T$ for
some invertible matrix $P$.  The main classical theorem on this topic
says that if 
$\chara(F)\neq 2$ then every symmetric bilinear form can
be diagonalized by a change of basis.  The question remains to decide
when two given diagonal matrices $D_1$ and $D_2$ represent equivalent
bilinear forms.  For instance, do $\begin{bmatrix} 2 & 0 \\ 0 & 11
\end{bmatrix}$ and $\begin{bmatrix} 3 & 0 \\ 0 & 1 \end{bmatrix}$
represent the same form over $\Q$?  

To pursue this question one looks for invariants.  The most obvious of
these is the rank of the matrix $A$.  This is in fact the unique
invariant when the field is algebraically closed.  For suppose a form
is represented by a diagonal matrix $D$, and let $\lambda$ be a
nonzero scalar. Construct a new basis by replacing the $i$th basis
element $e_i$ by $\lambda e_i$.  The matrix of the form with respect
to this new basis is the same as $D$, but with the $i$th diagonal
entry multiplied by $\lambda^2$.  The conclusion is that multiplying
the entries of $D$ by squares does not change the isomorphism class of
the underlying form.  This leads immediately to the classical theorem
saying that if every element of $F$ is a square (which we'll write as
$F=F^2$) then a symmetric bilinear form is completely classified by
its rank.

We now restrict to nondegenerate forms, in which case the matrix $A$
is nonsingular.  The element $\det(A)\in F^*$ is not quite an
invariant of the bilinear form, since after a change of basis the
determinant of the new matrix will be
$\det(P)\det(A)\det(P^T)=\det(P)^2\det(A)$.  However, the determinant
is a well-defined invariant if we regard it as an element of
$F^*/(F^*)^2$.  Since $\frac{22}{3}$ is not a square in $\Q$, for
instance, this tells us that the matrices $\begin{bmatrix} 2 & 0 \\ 0
& 11 \end{bmatrix}$ and $\begin{bmatrix} 3 & 0 \\ 0 & 1 \end{bmatrix}$
don't represent isomorphic forms over $\Q$.

The rank and determinant are by far the simplest invariants to write
down, but they are not very strong.  They don't even suffice to
distinguish forms over $\R$.  This case is actually a good example to
look at.  For $a_1,\ldots,a_n \in \R^*$, let $\br{a_1,\ldots,a_n}$
denote the form on $\R^n$ defined by $(e_i,e_j)=\delta_{i,j}a_i$.
Since every element of $\R$ is either a square or the negative of a
square, it follows that every nondegenerate real form is isomorphic to
an $\br{a_1,\ldots,a_n}$ where each $a_i\in \{1,-1\}$.  When are two
such forms isomorphic?  Of course one knows the answer, but let's
think through it.  The Witt Cancellation Theorem (true over any field)
says that if $\br{x_1,\ldots,x_n,y_1,\ldots,y_k} \iso
\br{x_1,\ldots,x_n,z_1,\ldots,z_k}$ then $\br{y_1,\ldots,y_k} \iso
\br{z_1,\ldots,z_k}$.  So our problem reduces to deciding whether the
$n$-dimensional forms $\br{1,1,\ldots,1}$ and $\br{-1,\ldots,-1}$ are
isomorphic.  When $n$ is odd the determinant distinguishes them, but
when $n$ is even it doesn't.  Of course the thing to say is that the
associated quadratic form takes only positive values in the first
case, and only negative values in the second---but this is not exactly
an `algebraic' way of distinguishing the forms, in that it uses the
ordering on $\R$ in an essential way.  By the end of this section we
will indeed have purely algebraic invariants we can use here.

\subsection{The Grothendieck-Witt ring}
\label{se:GW}
In a moment we'll return to the problem of finding invariants more
sophisticated than the rank and determinant, but first we need a
little more machinery.  From now on $\chara(F)\neq 2$.  By a
\dfn{quadratic space} I mean a pair $(V,\mu)$ consisting of a
finite-dimensional vector space and a non-degenerate bilinear form
$\mu$.  To systemize their study one defines the {\it
Grothendieck-Witt ring\/} $GW(F)$.  This is the free abelian group
generated by isomorphism classes of pairs $(V,\mu)$, with the usual
relation identifying the direct sum of quadratic spaces with the sum
in the group.  The multiplication is given by tensor product of vector
spaces.

The classical theory of bilinear forms allows us to give a complete
description of the abelian group $GW(F)$ in terms of generators and
relations.  Recall that $\br{a_1,\ldots,a_n}$ denotes the
$n$-dimensional space $F^n$ with $( e_i,e_j)=\delta_{ij}a_i$.  So
$\br{a_1,\ldots,a_n}=\br{a_1}+\cdots+\br{a_n}$ in $GW(F)$.  The
fact that every symmetric bilinear form is diagonalizable tells us
that $GW(F)$ is generated by the elements $\br{a}$ for $a\in F^*$, and
we have already observed the relation $\br{ab^2}=\br{a}$ for any
$a,b\in F^*$.  As an easy exercise, one can also give a complete
description for when {\it two-dimensional\/} forms are isomorphic: one
must be able to pass from one to the other via the two relations
\begin{myequation}
\label{eq:GWrels}
 \br{ab^2}=\br{a} \quad \text{and}\quad \br{a,b}=\br{a+b,ab(a+b)}
\end{myequation}
where in the second we assume $a,b\in F^*$ and $a+b\neq 0$.
As an example, working over $\Q$ we have
\[ \br{3,-2}=\br{12,-2}=\br{10,-240}=\br{90,-15}.
\]
To completely determine all relations in $GW(F)$, one shows that if
two forms $\br{a_1,\ldots,a_n}$ and $\br{b_1,\ldots,b_n}$ are
isomorphic then there is a chain of isomorphic diagonal forms
connecting one to the other, where each link of the chain differs in
exactly two elements.  Thus, (\ref{eq:GWrels}) is a complete set of
relations for $GW(F)$.  The reader may consult \cite[2.9.4]{S1} for
complete details here.

The multiplication in $GW(F)$ can be described compactly by
\[ \br{a_1,\ldots,a_n}\cdot \br{b_1,\ldots,b_k}=\sum_{i,j}
\br{a_ib_j}.
\]

\subsection{The Witt ring}
\label{se:Witt}
The {\it Witt ring\/} $W(F)$ is the quotient of $GW(F)$ by the ideal
generated by the so-called `hyperbolic plane' $\br{1,-1}$.
Historically $W(F)$ was studied long before $GW(F)$, probably because
it can be defined without formally adjoining additive inverses as was
done for $GW(F)$.  One can check that the forms $\br{a,-a}$ and
$\br{1,-1}$ are isomorphic, and therefore if one regards hyperbolic
forms as being zero then $\br{a_1,\ldots,a_n}$ and
$\br{-a_1,\ldots,-a_n}$ are additive inverses.  So $W(F)$ can be
described as a set of equivalence classes of quadratic spaces,
and doesn't require working with `virtual' objects.

Because $\br{a,-a}\iso\br{1,-1}$ for any $a$, it follows that the ideal
$\bigl (\br{1,-1} \bigr )$ is precisely the additive subgroup of
$GW(F)$ generated by $\br{1,-1}$.  As an abelian group, it is just a
copy of $\Z$.  So we have the exact sequence $0 \ra \Z \ra GW(F) \ra
W(F) \ra 0$.

Let $GI(F)$ be the kernel of the dimension function $\dim\colon
GW(F)\ra \Z$, usually called the augmentation ideal.  Let $I(F)$ be
the image of the composite $GI(F)\inc GW(F) \fib W(F)$; one can check
that $I(F)$ consists precisely of equivalence classes of
even-dimensional quadratic spaces.  Note that $I$ is additively
generated by forms $\br{1,a}$, and therefore $I^n$ is additively
generated by $n$-fold products $\br{1,a_1}\br{1,a_2}\cdots
\br{1,a_n}$.  

The dimension function gives an isomorphism $W/I \ra \Z/2$.  The
determinant gives us a group homomorphism $GW(F) \ra F^*/(F^*)^2$, but
it does not extend to the Witt ring because $\det \br{1,-1}=-1$.  One
defines the {\it discriminant\/} of $\br{a_1,\ldots,a_n}$ to be
$(-1)^{\frac{n(n-1)}{2}}\cdot (a_1\cdots a_n)$, and with this
definition the discriminant gives a map of sets $W(F) \ra
F^*/(F^*)^2$.  It is not a homomorphism, but if we restrict to $I(F)
\ra F^*/(F^*)^2$ then it {\it is\/} a homomorphism.  As the
discriminant of $\br{1,a}\br{1,b}$ is a square, the elements of $I^2$
all map to $1$.  So we get an induced map $I/I^2 \ra F^*/(F^*)^2$,
which is obviously surjective.  It is actually an isomorphism---to see
this, note that 
\[ \br{x,y}\br{-1,y}=\br{-x,xy,-y,y^2}=\br{1,-x,-y,xy}
\]
and so $\br{x,y}\equiv\br{1,xy}$ (mod $I^2$).  It follows inductively
that $\br{a_1,\ldots,a_{2n}} \equiv \br{1,1,\ldots,1,a_1a_2\cdots
a_{2n}}$ (mod $I^2$).  So if $\br{a_1,\ldots,a_{2n}}$ is a form whose
discriminant is a square, it is equivalent mod $I^2$ to either
$\br{1,1,\ldots,1}=2n\br{1}$ (if $n$ is even) or
$\br{1,1,\ldots,1,-1}=(2n-2)\br{1}$ (if $n$ is odd).  In the former
case $2n\br{1}=2\br{1} \cdot n\br{1} \in I^2$, and in the latter case
$(2n-2)\br{1}=2\br{1} \cdot (n-1)\br{1} \in I^2$.  In either case we
have $\br{a_1,\ldots,a_{2n}}\in I^2$, and this proves injectivity.

The examples in the previous paragraph are very special, but
they suggest why one might hope for `higher' invariants which give
isomorphisms between the groups $I^n/I^{n+1}$ and something more
explicitly defined in terms of the field $F$.  This is what the Milnor
conjecture is about.

\begin{remark}
\label{re:GI}
For future reference, note that $2\br{1}=\br{1,1}\in I$, and therefore
the groups $I^n/I^{n+1}$ are $\Z/2$-vector spaces.  Also observe
that $GI(F)$ does not intersect the kernel of $GW(F)\ra W(F)$, and so
$GI(F)\ra I(F)$ is an isomorphism.  It follows that $(GI)^n/(GI)^{n+1}
\iso I^n/I^{n+1}$, for all $n$.
\end{remark}

\vspace{0.1in}

\subsection{More invariants}
Recall that the \dfn{Brauer group} $\Br(F)$ is a set of equivalence
classes of central, simple $F$-algebras, with the group structure
coming from tensor product.  The inverse of such an algebra is its
opposite algebra, where the order of multiplication has been reversed.

From a quadratic space $(V,\mu)$ one can construct the associated
Clifford algebra $C(\mu)$: this is the quotient of the tensor algebra
$T_F(V)$ by the relations generated by $v\tens v=\mu(v,v)$.  Clifford
algebras are $\Z/2$-graded by tensor length.  If $\mu$ is
even-dimensional then $C(\mu)$ is a central simple algebra, and if
$\mu$ is odd-dimensional then the even part $C_0(q)$ is a central
simple algebra.  So we get an invariant of quadratic spaces taking its
values in $\Br(F)$ (see \cite[9.2.12]{S1} for more detail).  This is
usually called the \dfn{Clifford invariant}, or sometimes the
\dfn{Witt invariant}.  Since any Clifford algebra is isomorphic to its
opposite, the invariant always produces a $2$-torsion class.

Now we need to recall some Galois cohomology.  Let $\bar{F}$ be a
separable closure of $F$, and let $G=\Gal(\bar{F}/F)$.  Consider the
short exact sequence of $G$-modules $0\ra \Z/2 \ra \bar{F}^* \ra \bar{F}^*
\ra 0$, where the second map is squaring.  Hilbert's Theorem
90 implies that $H^1(G;\bar{F}^*)=0$, which means that the induced
long exact sequence in Galois cohomology splits up into
\[ 0 \ra H^0(G;\Z/2) \ra {F}^* \llra{2} {F}^* \ra H^1(G;\Z/2)
\ra 0 \]
and
\[ 0 \ra H^2(G;\Z/2) \ra H^2(G;\bar{F}^*) \llra{2} H^2(G;\bar{F}^*).
\]
The group $H^2(G;\bar{F}^*)$ is known to be isomorphic to $\Br(F)$, so
we have $H^0(G;\Z/2)=\Z/2$, $H^1(G;\Z/2)=F^*/(F^*)^2$, and the
$2$-torsion in the Brauer group is precisely $H^2(G;\Z/2)$.  From now
on we will write $H^*(F;\Z/2)=H^*(G;\Z/2)$.

At this point we have the rank map $e_0\colon W(F) \ra
\Z/2=H^0(F;\Z/2)$, which gives an isomorphism $W/I \ra \Z/2$.  We have
the discriminant $e_1\colon I(F) \ra F^*/(F^*)^2=H^1(F;\Z/2)$ which
gives an isomorphism $I/I^2 \ra F^*/(F^*)^2$, and we have the Clifford
invariant $e_2\colon I^2 \ra H^2(F;\Z/2)$.  With a little work one can
check that $e_2$ is a homomorphism, and it kills $I^3$.  The question
of whether $I^2/I^3 \ra H^2(F;\Z/2)$ is an isomorphism is difficult,
and wasn't proven until the early 80s by Merkurjev \cite{M} (neither
surjectivity nor injectivity is obvious).  The maps $e_0,e_1,e_2$ are
usually called the {\it classical invariants\/} of quadratic forms.

The above isomorphisms can be rephrased as follows. The ideal $I$
consists of all elements where $e_0=0$; $I^2$ consists of all elements
such that $e_0=0$ and $e_1=1$; and by Merkujev's theorem $I^3$ is
precisely the set of elements for which $e_0$, $e_1$, and $e_2$ are all
trivial.  Quadratic forms will be completely classified by these
invariants if $I^3=0$, but unfortunately this is usually not the case.
This brings us to the search for higher invariants.  One early result
along these lines is due to Delzant \cite{De}, who defined
Stiefel-Whitney invariants with values in Galois
cohomology.  Unfortunately these are not the `right' invariants, as
they do not lead to complete classifications for elements in $I^n$,
$n\geq 3$.  

\subsection{Milnor's work}
At this point we find ourselves looking at the two rings $\Gr_I W(F)$
and $H^*(F;\Z/2)$, and we have maps between them in dimensions $0$,
$1$, and $2$.  I think Milnor, inspired by his work on algebraic
$K$-theory, wrote down the best ring he could find which would map to
both rings above.  In \cite{Mi} he defined what is now called `Milnor
$K$-theory' as
\[ K^M_*(F)=T_\Z(F^*)/\lb a\tens (1-a)| a\in F-\{0,1\}\rb 
\]
where $T_\Z(V)$ denotes the tensor algebra over $\Z$ on the abelian
group $V$.  The grading comes from the grading on the tensor algebra,
in terms of word length.  I will write $\brn{a_1,\ldots,a_n}$ for the
element $a_1 \tens \cdots \tens a_n \in K^M_n(F)$.  

Note that when dealing with $K^M_*(F)$ one must be careful not to
confuse the addition---which comes from  multiplication in
$F^*$---with the multiplication.  So for instance
$\brn{a}+\brn{b}=\brn{ab}$ but $\brn{a}\cdot\brn{b}=\brn{a,b}$. 
This is in contrast to the operations in $GW(F)$, where one has
$\br{a}+\br{b}=\br{a,b}$ and $\br{a}\tens\br{b}=\br{ab}$.
Unfortunately it is very easy to get these confused. 
Note that $\brn{a^2}=2\brn{a}$, and more generally
$\brn{a^2,b_1,\ldots,b_n}=2\brn{a,b_1,\ldots,b_n}$.  

\begin{remark}
From a modern perspective the name `K-theory' applied to $K^M_*(F)$ is
somewhat of a misnomer; one should not take it too seriously.  The
construction turns out to be more closely tied to algebraic cycles than to
algebraic $K$-theory, and so I personally like the term `Milnor cycle
groups'.  I doubt this will ever catch on, however.
\end{remark}

Milnor produced two ring homomorphisms $\eta\colon K^M_*(F)/2 \ra
H^*(F;\Z/2)$ and $\nu \colon K^M_*(F)/2 \ra \Gr_I W(F)$.  To define
the map $\nu$, note first that we have already established an
isomorphism $F^*/(F^*)^2 \ra I/I^2$ sending $\brn{a}$ to
$\br{a,-1}=\br{a}-\br{1}$ (this is the inverse of the discriminant).
This tells us what $\nu$ does to elements in degree $1$.  Since these
elements generate $K^M_*(F)$ multiplicatively, to construct $\nu$ it
suffices to verify that the appropriate relations are satisfied in the
image.  So we first need to check that
\[ 0=\Bigl (\br{a}-\br{1} \Bigr )\cdot \Bigl (\br{1-a}-\br{1} \Bigr) =
\br{a(1-a)}-\br{a}-\br{1-a}+\br{1} = \br{a(1-a),1}-\br{a,1-a},
\]
but this follows directly from the second relation in
(\ref{eq:GWrels}).  We also must check that $2\brn{a}$ maps to $0$,
but $2\brn{a}=\brn{a^2} \mapsto \br{a^2}-\br{1}$ and the latter
vanishes by the first relation in (\ref{eq:GWrels}).  For future
reference, note that $\nu(\brn{a})$ is equal to both $\br{a,-1}$ and
$\br{-a,1}$ in $I/I^2$, since this group is $2$-torsion.

Defining $\eta$ is similar.  We have already noticed that
there is a natural isomorphism $H^1(F;\Z/2)\iso F^*/(F^*)^2$, and so
it is clear where the element $\brn{a}$ in $K^M_1(F)=F^*$ must be sent.  The
verification that $a \cup (1-a)=0$ in $H^2(F;\Z/2)$ is in
\cite[6.1]{Mi}.

Milnor observed that both $\eta$ and $\nu$ were isomorphisms in all
the cases he could compute.  The claim that $\eta$ is an isomorphism
is nowadays known as {\it the Milnor conjecture\/}, and was proven by
Voevodsky in 1996 \cite{Vmilnor}.  The claim that $\nu$ is an
isomorphism goes under the name {\it Milnor's conjecture on quadratic
forms\/}.  For characteristic $0$ it was proven in 1996 by Orlov,
Vishik, and Voevodsky \cite{OVV}, who deduced it as a consequence of
the work in \cite{Vmilnor}.  I believe the proof now works in
characteristic $p$, based on the improved results of \cite{Vmilnor2}.
A second proof, also in characteristic $0$, was outlined by Morel
\cite{M2} using the motivic Adams spectral sequence, and again
depended on results from \cite{Vmilnor}; unfortunately complete
details of Morel's proof have yet to appear.

It is interesting that the conjecture on quadratic forms doesn't have
an independent proof, and is the less primary of the two.  Note that
both $K^M_*(F)/2$ and $GW(F)$ can be completely described in terms of
generators and relations (although the latter does not quite imply
that we know all the relations in $\Gr_I W(F)$, which is largely the
problem).  The map $\nu$ is easily seen to be surjective, and so the
only question is injectivity.  Given this, it is in some ways
surprising that the conjecture is as hard as it is.

\begin{remark}
The map $\eta$ is called the {\it norm residue symbol\/}, and can be
defined for primes other than $2$.  The {\it Bloch-Kato conjecture\/}
is the statement that $\eta\colon K^M_i(F)/l \ra H^i(F;\mu_l^{\tens
i})$ is an isomorphism for $l$ a prime different from $\chara(F)$.
This is a direct generalization of the Milnor conjecture to the case
of odd primes.  A proof was released by Voevodsky in 2003
\cite{Vbkato} (although certain auxiliary results required for the
proof remain unwritten).  I'm not sure anyone has ever considered an
odd-primary analog of Milnor's conjecture on quadratic forms---what
could replace the Grothendieck-Witt ring here?
\end{remark}

At this point it might be useful to think through the Milnor
conjectures in a few concrete examples.  For these we refer the reader
to Appendix A.  Let's at least note here that through the work of
Milnor, Bass, and Tate (cf. \cite{Mi}) the conjectures could be
verified for all finite fields and for all finite extensions of
$\Q$ (in fact for all global fields).

Finally, let's briefly return to the classification of forms over
$\R$.  We saw earlier that this reduces to proving that the
$n$-dimensional forms $\br{1,1,\ldots,1}$ and $\br{-1,-1,\ldots,-1}$
are not isomorphic.  Can we now do this algebraically?  If they were
isomorphic, they would represent the same element of $W(\R)$.  It
would follow that $(2n)\br{1}=0$ in $W(\R)$.  Can this happen?  The
isomorphisms $\Z/2[a] \iso H^*(\Z/2;\Z/2) \iso K^M_*(\R)/2 \iso \Gr_I
W(\R)$ show that $\Gr_I W(\R)$ is a polynomial algebra on the class
$\br{-1,-1}$ (the generator $a$ corresponds to the generator $-1$ of
$\R^*/(\R^*)^2$, and $\nu(-1)=\br{-1,-1}$).  It follows that
$2^k\br{1}=\pm \br{-1,-1}^k$ is a generator for the group
$I^k/I^{k+1}\iso \Z/2$.  If $m=2^ir$ where $r$ is odd, then
$m\br{1}=2^i\br{1} \cdot r\br{1}$.  Since $r\br{1}$ is the generator
for $W/I$ and $2^i\br{1}$ is a generator for $I^i/I^{i+1}$, it follows
that $m\br{1}$ is also a generator for $I^i/I^{i+1}$.  In particular,
$m\br{1}$ is nonzero.  So we have proven via algebraic methods
(although in this case also somewhat pathological ones) that
$\br{1,1,\ldots,1}\not\iso \br{-1,-1,\ldots,-1}$.

\subsection{Further background reading}
There are several good expository papers on the theory of quadratic
forms, for example \cite{Pf1} and \cite{S2}.  The book \cite{S1} is a
very thorough and readable resource as well.  For the Milnor
conjectures themselves there is \cite{Pf2}, which in
particular gives several applications of the conjectures; it
also gives detailed references to original papers.  The
beginning sections of \cite{AEJ} offer a nice survey concerning the
search for `higher' invariants of quadratic forms.  It's worth pointing
out that after Milnor's work definitions of $e_3$, $e_4$, and $e_5$
were eventually given---with a lot of hard work---but this was the
state of the art until 1996.  Finally, the introduction of
\cite{Vmilnor2} gives a history of work on the Milnor conjecture.

\vfill\eject


\section{Proof of the conjecture on the norm residue symbol}
This section outlines Voevodsky's proof of the Milnor conjecture on
the norm residue symbol \cite{Vmilnor,Vmilnor2}.  Detailed,
step-by-step summaries have been given in \cite{M1} and \cite{Su}.  My
intention here is not to give a complete, mathematically rigorous
presentation, but rather just to give the flavor of what is involved.

Several steps in the proof involve manipulations with motivic
cohomology based on techniques that were developed in \cite{VSF}.  I
have avoided giving any details about these steps, in an attempt to
help the exposition.  Most of these details are not hard to
understand, however---there are only a few basic techniques to keep
track of, and one can read about them in \cite{VSF} or \cite{MVW}.
But I hope that by keeping some of this stuff in a black box the
overall structure of the argument will become clearer.

\subsection{Initial observations}
The aim is to show that $\eta\colon K^M_*(F)/2 \ra H^*(F;\Z/2)$ is an
isomorphism.  To do this, one of the first things one might try to
figure out is what kind of extra structure $K^M_*(F)/2$ and
$H^*(F;\Z/2)$ have in common.  For instance, they are both covariant
functors in $F$, and the covariance is compatible with the norm
residue symbol.  It turns out they both have transfer maps for finite
separable extensions (which, for those who like to think
geometrically, are the analogs of covering spaces).  That is, if
$j\colon F \inc F'$ is a separable extension of degree $n$ then there
is a map $j^{!}\colon K^M_*(F') \ra K^M_*(F)$ such that $j^!j_*$ is
multiplication by $n$, and similarly for $H^*(F;\Z/2)$.  (Note that
the construction of transfer maps for Milnor $K$-theory is not at all
trivial---some ideas were given in \cite[Sec. 5.9]{BT}, but the full
construction is due to Kato \cite[Sec. 1.7]{Ka}).  It follows that if
$n$ is odd then $K^M_*(F)/2 \ra H^*(F;\Z/2)$ is a retract of the map
$K^M_*(F')/2 \ra H^*(F';\Z/2)$.  So if one had a counterexample to the
Milnor conjecture, field extensions of all odd degrees would still be
counterexamples.  This is often referred to as ``the transfer
argument''.

Another observation is that both functors can be extended to rings
other than fields, and if $R$ is a discrete valuation ring then both
functors have a `localization sequence' relating their values on $R$,
the residue field, and the quotient field.  I will not go into details
here, but if $F$ is a field of characteristic $p$ then by using the
Witt vectors over $F$ and the corresponding localization sequence, one
can reduce the Milnor conjecture to the case of characteristic $0$
fields.  The argument is in \cite[Lemma 5.2]{Vmilnor}.  In Voevodsky's
updated proof of the Milnor conjecture \cite{Vmilnor2} this step is
not necessary, but I think it's useful to realize that the Milnor
conjecture is not hard because of `crazy' things that might happen in
characteristic $p$---it is hard even in
characteristic $0$.

\subsection{A first look at the proof}
\label{se:firstlook}
The proof goes by induction.  We assume the norm residue map
$\eta\colon K_*^M(F)/2 \ra H^*(F;\Z/2)$ is an isomorphism for all
fields $F$ and all $*<n$,
and then prove it is also an isomorphism for $*=n$.  The basic theme
of the proof, which goes back to Merkurjev, involves two steps:
\begin{enumerate}[(1)]
\item Verify that $\eta_n$ is an isomorphism for certain `big'
fields---in our case, those which have no extensions of odd degree and
also satisfy $K_n(F)=2K_n(F)$ (so that one must prove
$H^n(F;\Z/2)=0$).  Notice that when $n=1$ the condition $K_1=2K_1$
says that $F=F^2$.

\item Prove that if $F$ were a field for which $\eta_n$ is {\it not\/}
an isomorphism then one could expand $F$ to make a `bigger'
counterexample, and could keep doing this until you're in the range
covered by step (1).  This would show that no such $F$ could exist.

In more detail one shows that for any $\brn{a_1,\ldots,a_n}\in K_n(F)$
one can construct an extension $F\inc F'$ with the property that
$\brn{a_1,\ldots,a_n}\in 2K_n(F')$ and $\eta_n\colon K_n(F')/2 \ra
H^n(F';\Z/2)$ still fails to be an isomorphism.  By doing this over
and over and taking a big colimit, one gets a counterexample
where $K^M_n=2K^M_n$.  
\end{enumerate} 

\vspace{0.1in}

Neither of the above two steps is trivial, but step (1) involves
nothing very fancy---it is a calculation in Galois cohomology
which takes a few pages, but is not especially hard.  See
\cite[Section 5]{Vmilnor2}.  Step (2) is the more subtle and
interesting step.  Note that if $\und{a}=\brn{a_1,\ldots,a_n} \notin
2K^M_n(F)$ then none of the $a_i$'s can be in $F^2$.  There are
several ways one can extend $F$ to a field $F'$ such that $\und{a} \in
2K^M_n(F')$: one can adjoin a square root of any $a_i$, for
instance.  The problem is to find such an extension where you have
enough control over the horizontal maps in the diagram
\[ \xymatrix{
 K^M_n(F)/2 \ar[r]\ar[d]_{\eta_F} & K^M_n(F')/2 \ar[d]^{\eta_{F'}} \\
 H^n(F;\Z/2) \ar[r] & H^n(F';\Z/2)
}
\]
to show that if $\eta_F$ fails to be an isomorphism then so does
$\eta_{F'}$.  The selection of the `right' $F'$ is delicate.

We will alter our language at this point, because we will want to
bring more geometry into the picture.  Any finitely-generated
separable extension $F\inc F'$ is the function field of a smooth
$F$-variety.  A \dfn{splitting variety} for an element $\und{a}\in
K^M_n(F)$ is a smooth variety $X$, of finite type over $F$, with the
property that $\und{a} \in 2K^M_n(F(X))$.  Here $F(X)$ denotes the
function field of $X$.  As we just remarked, there are many such
varieties: $X=\Spec F[u]/(u^2-a_1)$ is an example.  The particular
choice we'll be interested in is more complicated.

Given $b_1,\ldots,b_k\in F$, let $q_{\und{b}}$ be the quadratic form
in $2^k$ variables corresponding to the element
\[ \br{1,-b_1}\tens \br{1,-b_2} \tens \cdots \tens \br{1,-b_k} \in GW(F).
\]
For example,
$q_{b_1,b_2}(x_1,\ldots,x_4)=x_1^2-b_1x_2^2-b_2x_3^2+b_1b_2 x_4^2$.
Such $q$'s are called \dfn{Pfister forms}, and they have a central
role in the modern theory of quadratic forms (see \cite[Chapter
4]{S1}, for instance).

For $a_1,\ldots,a_n\in F$, define $Q_{\und{a}}$ to be the projective
quadric in $\P^{2^{n-1}}$ given by the equation
\[ q_{a_1,\ldots,a_{n-1}}(x_0,\ldots,x_{[2^{n-1}-1]})-a_nx_{2^{n-1}}^2=0.
\]
In \cite{Vmilnor2} these are called
\dfn{norm quadrics}.  A routine argument
\cite[Prop. 4.1]{Vmilnor2} shows that $Q_{\und{a}}$ is a splitting variety
for $\und{a}$.  The reason for choosing to study this particular
splitting variety will not be clear until later; isolating this object
is one of the key aspects of the proof.

The name of the game will be to understand enough about the difference
between $K_n^M(F)/2$ and $K_n^M(F(Q_{\und a}))/2$ (as well as the
corresponding Galois cohomology groups) to show that $K^M_n(F(Q_{\und
a}))/2 \ra H^n(F(Q_{\und a});\Z/2)$ still fails to be an isomorphism.
Voevodsky's argument uses motivic cohomology---of the quadrics
$Q_{\und a}$ and other objects---to `bridge the gap' between
$K_n^M(F)/2$ and $K^M_n(F(Q_{\und a}))/2$.  

\subsection{Motivic cohomology enters the picture}
Motivic cohomology is a bi-graded functor $X\mapsto H^{p,q}(X;\Z)$
defined on the category of smooth $F$-schemes.  Actually it is defined
for all simplicial smooth schemes, as well as for more general
objects.  One of the lessons of the last ten years is that one can set
up a model category which contains all these objects, and then a
homotopy theorist can deal with them in much the same ways he deals
with ordinary topological spaces.  From now on I will do this
implicitly (without ever referring to the machinery involved).

The coefficient groups $H^{p,q}(\Spec F;\Z)$ vanish for $q<0$ and for
$p>q\geq 0$.  For us an important point is that the groups
$H^{n,n}(\Spec F;\Z)$ are canonically isomorphic to $K^M_n(F)$.
Proving this is not simple!  An account is given in \cite[Lecture
5]{MVW}.  Finally, we note that one can talk about motivic
cohomology with finite coefficients $H^{p,q}(X;\Z/n)$, related to
integral cohomology via the exact sequence
\[ \cdots \ra 
H^{p,q}(X;\Z) \llra{\times n} H^{p,q}(X;\Z) \ra H^{p,q}(X;\Z/n) \ra
H^{p+1,q}(X;\Z) \ra \cdots 
\]
The sequence shows  $H^{n,n}(\Spec F;\Z/2)\iso K^M_n(F)/2$ and
$H^{p,q}(\Spec F;\Z/2)=0$ for $p> q \geq 0$.

Now, there is also an analagous theory $H^{p,q}_L(X;\Z)$ which is
called \dfn{Lichtenbaum (or \'etale) motivic cohomology}.  There is a
natural transformation $H^{p,q}(X;\Z) \ra H^{p,q}_L(X;\Z)$.  The
theory $H^{*,*}_L$ is the closest theory to $H^{*,*}$ which satisfies
descent for the \'etale topology (essentially meaning that when $E\ra
B$ is an \'etale map there is a spectral sequence starting with
$H^{*,*}_L(E)$ and converging to $H^{*,*}_L(B)$).  The relation
between $H^{*,*}$ and $H^{*,*}_L$ is formally analagous to that
between a cohomology theory and a certain Bousfield localization of
it.  It is known that $H^{p,q}_L(X;\Z/n)$ is canonically isomorphic to
\'etale cohomology $H^{p}_{et}(X;\mu_n^{\tens q})$, if $n$ is prime to
$\chara(F)$.  From this it follows that $H_L^{p,q}(\Spec F;\Z/2)$ is the
Galois cohomology group $H^p(F;\Z/2)$, for all $q$.  At this point we
can re-phrase the Milnor conjecture as the statement that the maps
$H^{p,p}(\Spec F;\Z/2) \ra H^{p,p}_{L}(\Spec F;\Z/2)$ are
isomorphisms.

There are other conjectures about the relation between $H^{*,*}$
and $H_L^{*,*}$ as well.  A conjecture of Lichtenbaum says that
$H^{p,q}(X;\Z) \ra H^{p,q}_L(X;\Z)$ should be an isomorphism
whenever $p\leq q+1$.  Note that this would imply a corresponding
statement for $\Z/n$-coefficients, and in particular would imply the
Milnor conjecture.  Also, since one knows $H^{n+1,n}(\Spec F;\Z)=0$ 
Lichtenbaum's conjecture  would
imply that $H^{n+1,n}_L(\Spec F;\Z)$ also vanishes.  This latter
statement was conjectured independently by both Beilinson and
Lichtenbaum, and is known as a the \dfn{Generalized Hilbert's Theorem
90} (the case $n=1$ is a translation of the statement that
$H^1_{Gal}(F;\bar{F}^*)=0$, which follows from the classical Hilbert's
Theorem 90).

By knowing enough about how to work with motivic cohomology, Voevodsky
was able to prove the following relation among these conjectures 
(as well as other relations which we won't need):

\begin{prop}
\label{pr:reduce}
Fix an $n\geq 0$.  Assume that $H^{k+1,k}_L(\Spec F;\Zlt)=0$ for all
fields $F$ and all $0\leq k \leq n$. Then for any smooth simplicial
scheme $X$ over a field $F$, the maps $H^{p,q}(X;\Z/2) \ra
H^{p,q}_L(X;\Z/2)$ are isomorphisms when $q\geq 0$ and $p\leq q \leq
n$; and they are monomorphisms for $p-1=q\leq n$.  In particular,
applying this when $p=q$ and $X=\Spec F$ verifies the Milnor conjecture
in dimensions $\leq n$.
\end{prop}

It's worth pointing out that the proof uses nothing special about
the prime $2$, and so the statement is valid for all other primes as
well.  

For us, the importance of the above proposition is two-fold.  First,
it says that to prove the Milnor conjecture one only has to worry
about the vanishing of one set of groups (the $H^{n+1,n}_L$'s) rather
than two sets (the kernel and cokernel of $\eta$).  Secondly,
inductively assuming that the Generalized Hilbert's Theorem 90 holds
up through dimension $n$ is going to give us a lot more to work with
than inductively assuming the Milnor conjecture up through dimension
$n$.  Instead of just knowing stuff about $H^{n,n}$ of fields, we know
stuff about $H^{p,q}$ of any smooth simplicial scheme.  The need for
this extra information is a key feature of the proof.

\vspace{0.1in}

\subsection{\CCech complexes}
We only need one more piece of machinery before returning to the proof
of the Milnor conjecture.  This piece is hard to motivate, and its
introduction is one of the more ingenious aspects of the proof.  The
reader will just have to wait and see how it arises in
section~\ref{se:proof} (see also Remark~\ref{re:cech}).

Let $X$ be any scheme.  The \dfn{\CCech complex} \mdfn{$\Cech X$} is
the simplicial scheme with $(\Cech X)_n=X \times X \times \cdots
\times X$ ($n+1$ factors) and the obvious face and degeneracies.
This simplicial scheme can be regarded as augmented by the map $X \ra
\Spec F$.  

For a topological space the realization of the associated \CCech
complex is always contractible---in fact, choosing any point of
$X$ allows one to write down a contracting homotopy for the simplicial
space $\Cech X$.  But in algebraic geometry the scheme $X$ may not
have rational points; i.e., there may not exist any maps $\Spec F \ra
X$ at all!  If $X$ does have a rational point then the same trick lets
one write down a contracting homotopy, and therefore $\Cech X$ behaves
as if it were $\Spec F$ in all computations.  (More formally, $\Cech X$ is
homotopy equivalent to $\Spec F$ in the motivic homotopy category).

Working in the motivic homotopy category, one finds that for any
smooth scheme $Y$ the set of homotopy classes $[Y,\Cech X]$ is either
empty or a singleton.  The latter holds precisely if $Y$ admits a
Zariski cover $\{U_\alpha\}$ such that there exist scheme maps
$U_\alpha \ra X$ (not necessarily compatible on the
intersections).  The object $\Cech X$ has no `higher homotopy
information', only this very simple discrete information about whether
or not certain maps exist.  One should think of $\Cech X$ as very
close to being contractible.  I point out again that in topology there
is always at least one map between nonempty spaces, and so $\Cech X$
is not very interesting.

If $E\ra B$ is an \'etale cover, then there is a spectral sequence
whose input is $H^{*,*}_L(E;\Z)$ and which converges to
$H^{*,*}_L(B;\Z)$ (this is the \'etale descent property).  In
particular, if $X$ is a smooth scheme and we let $F'=F(X)$,
$X'=X\times_F F'$, then $X' \ra X$ and $\Spec F'\ra \Spec F$ are both
\'etale covers.  The scheme $X'$ necessarily has a rational point over
$F'$, so  $\Cech X'$
and $\Spec F'$ look the same to $H_L$.  The \'etale descent property
then shows that $\Cech X$ and $\Spec F$ also look the same: in other
words, the maps $H^{p,q}_L(\Spec F;\Z) \ra H^{p,q}_L(\Cech X;\Z)$ are
all isomorphisms (and the same for finite coefficients).  This is {\it
not\/} true for $H^{*,*}$ in place of $H^{*,*}_L$.  
One might paraphrase all this by saying that in the \'etale world $\Cech
X$ is contractible, just as it is in topology.

\vspace{0.1in}

\subsection{The proof}
\label{se:proof}
Now I am going to give a complete summary of the proof as it appears
in \cite{Vmilnor,Vmilnor2}.  Instead of proving the Milnor conjecture
in its original form one instead concentrates on the more manageable
conjecture that $H^{i+1,i}_L(\Spec F;\Zlt)=0$ for all $i$ and all
fields $F$.  One assumes this has been proven in the range $0\leq i
<n$, and then shows that it also follows for $i=n$.

Suppose that $F$ is a field with $H^{n+1,n}_L(F;\Zlt)\neq 0$.  The
transfer argument shows that any extension field of odd degree would
still be a counterexample, so we can assume $F$ has no extensions of
odd degree.  One checks via some Galois cohomology computations---see
\cite[section 5]{Vmilnor2}---that if such a field has $K_n^M(F)=2K_n^M(F)$
then $H^{n+1,n}_L(\Spec F;\Zlt)=0$.  So our counterexample cannot have
$K^M_n(F)=2K^M_n(F)$.  By the reasoning from
section~\ref{se:firstlook}, it will suffice to show that for every
$a_1,\ldots,a_n\in F$ the field $F(\Qa)$ is still a counterexample.
We will in fact show that $H^{n+1,n}_L(F;\Zlt) \ra
H^{n+1,n}_L(F(\Qa);\Zlt)$ is injective.

\vspace{0.1in}

Suppose $u$ is in the kernel of the above map, and consider the
diagram
\[ \xymatrix{
 & H^{n+1,n}_L(\Spec F;\Zlt) \ar[r]\ar[d]^\iso 
     & H^{n+1,n}_L(\Spec F(\Qa);\Zlt) \\
 H^{n+1,n}(\Cech \Qa;\Zlt) \ar[r] & H^{n+1,n}_L(\Cech \Qa;\Zlt).
}
\]
Let $u'$ denote the image of $u$ in $H^{n+1,n}_L(\Cech
\Qa;\Zlt)$. One can show (after some extensive manipulations with
motivic cohomology) that the hypothesis on $u$ implies that $u'$ is
the image of an element in $H^{n+1,n}(\Cech \Qa;\Zlt)$.  It will
therefore be sufficient to show that this group is zero.

Let $\CQa$ be defined by the cofiber sequence $(\Cech \Qa)_+ \ra
(\Spec F)_+ \ra \CQa$.  This means $\tH^{*,*}(\CQa)$ fits in an
exact sequence
\[ \ra H^{p-1,q}(\Cech \Qa) \ra \tH^{p,q}(\CQa) \ra H^{p,q}(\Spec F) 
\ra H^{p,q}(\Cech \Qa) \ra 
\tH^{p+1,q}(\CQa) \ra \cdots
\]
So the reduced motivic cohomology of $\CQa$ detects
the `difference' between the motivic cohomology of $\Cech \Qa$ and $\Spec
F$.  The fact that $H^{i,n}(\Spec F;\Z)=0$ for $i>n$ shows that
$H^{n+1,n}(\Cech \Qa;\Zlt)\iso \tH^{n+2,n}(\CQa;\Zlt)$.  Since $\Qa$
has a rational point (and therefore $\Cech \Qa$ is contractible) over
a degree $2$ extension of $F$, it follows from the transfer argument
that the above group is killed by $2$.  To show that the group is zero it
is therefore sufficient to prove that the image of
$\tH^{n+2,n}(\CQa;\Zlt) \ra \tH^{n+2,n}(\CQa;\Z/2)$ is zero.
This is the same as the image of
$\tH^{n+2,n}(\CQa;\Z) \ra \tH^{n+2,n}(\CQa;\Z/2)$, which 
I'll denote by $\intH^{n+2,n}(\CQa;\Z/2)$.

\vspace{0.1in}

So far most of what we have done is formal; but now we come to the
crux of the argument.  For any smooth scheme $X$ one has cohomology
operations acting on $H^{*,*}(X;\Z/2)$.  In particular, one can
produce analogs of the Steenrod operations: the Bockstein acts with
bi-degree $(1,0)$, and $\Sq^{2^i}$ acts with bi-degree
$(2^i,2^{i-1})$.  From these one defines the Milnor $Q_i$'s, which
have bi-degree $(2^{i+1}-1,2^i-1)$.  In ordinary topology these are
defined inductively by $Q_0=\beta$ and $Q_i=[Q_{i-1},\Sq^{2^i}]$,
whereas motivically one has to add some extra terms to this commutator
(these arise because the motivic cohomology of a point is nontrivial).
One shows that $Q_i\circ Q_i=0$, and that $Q_i=\beta q+q\beta$ for a
certain operation $q$.  It follows from the latter formula that $Q_i$
maps elements in $\intH$ to elements in $\intH$.  All of these facts
also work in ordinary topology, it's just that the proofs here are a
little more complex.
  
The next result is \cite[Cor. 3.8]{Vmilnor2}.  It is the first of two
main ingredients needed to complete the proof.

\begin{prop}
\label{pr:margolis}
Let $X$ be a smooth quadric in $\P^{2^n}$, and let $\tilde{C}X$ be
defined by the cofiber sequence $(\Cech X)_+ \ra (\Spec F)_+ \ra
\tilde{C}X$.  Then for $i\leq n$, every element of
$\tH^{*,*}(\tilde{C}X;\Z/2)$ that is killed by $Q_i$ is also in the
image of $Q_i$.
\end{prop}

This is a purely `topological' result, in that its proof uses no algebraic
geometry.  It follows from the most basic properties of the Steenrod
operations, motivic cohomology (like Thom isomorphism), and elementary
facts about the characteristic numbers of quadrics.  The argument is
purely homotopy-theoretic.

The second main result we will need is where all the algebraic
geometry enters the picture.  Voevodsky deduces it from results of
Rost, who showed that the motive of $\Qa$ splits off a certain direct
summand.  See \cite[Th. 4.9]{Vmilnor2}.

\begin{prop}
\label{pr:rost}
$\tH^{2^n,2^{n-1}}(\CQa;\Zlt)=0$.
\end{prop}

Using the above two propositions we can complete the proof of the
Milnor conjecture.  In order to draw a concrete picture, let us just
assume $n=4$ for the moment.  We are trying to show that
$\intH^{6,4}(\CQa;\Z/2)=0$.  Consider the
diagram
\[ \xymatrix{H^{p,q}(\Spec F;\Z/2) \ar[r]\ar[d] & H^{p,q}(\Cech
\Qa;\Z/2) \ar[d] \\ H^{p,q}_L(\Spec F;\Z/2) \ar[r]^\iso &
H^{p,q}_L(\Cech \Qa;\Z/2).  }\] Our inductive assumption together with
Proposition~\ref{pr:reduce} implies that the vertical maps are
isomorphisms for $p\leq q \leq n-1$, and monomorphisms for $p-1=q\leq
n-1$.  So the top horizontal map is an isomorphism in the first range
and a monomorphism in the second.  The long exact sequence in motivic
cohomology then shows that $\tH^{p,q}(\CQa;\Z/2)=0$ for $p\leq q \leq
n-1$.  This is where our induction hypothesis has gotten us.  The
following diagram depicts what we now know about
$\tH^{p,q}(\CQa;\Z/2)$ (the group marked $??$ is $\tH^{6,4}$, the one
we care about):

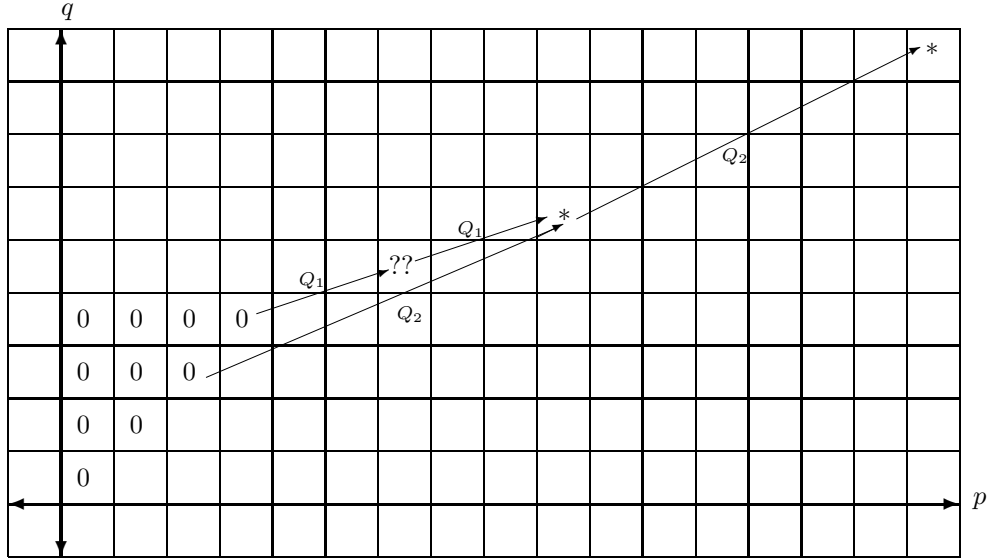
\begin{picture}(340,230)(12,-2)
\multiput(0,0)(20,0){19}{\line(0,1){200}}
\multiput(0,0)(0,20){11}{\line(1,0){360}}
\thicklines
%
\put(20,0){\vector(0,1){200}}
\put(20,200){\vector(0,-1){200}}
%
\put(0,20){\vector(1,0){360}}
\put(360,20){\vector(-1,0){360}}
\put(20,205){$q$}
\put(365,20){$p$}
\thinlines
\multiput(26,27)(20,0){1}{$0$}
\multiput(26,47)(20,0){2}{$0$}
\multiput(26,67)(20,0){3}{$0$}
\multiput(26,87)(20,0){4}{$0$}

\put(144,107){$??$}

\put(154,112){\vector(3,1){50}}
\put(94,92){\vector(3,1){50}}
\put(208,127){$*$}
\put(170,122){$\scriptstyle{Q_1}$}
\put(110,103){$\scriptstyle{Q_1}$}

\put(215,128){\vector(2,1){130}}
\put(347,190){$*$}
\put(270,150){$\scriptstyle{Q_2}$}
\put(75,68){\line(7,3){133}}
\put(200,121){\vector(2,1){10}}
\put(147,90){$\scriptstyle{Q_2}$}
\end{picture}

At this point Proposition~\ref{pr:margolis} shows that $Q_1\colon H^{6,4}
\ra H^{9,5}$ is injective, and that $Q_2\colon H^{9,5} \ra H^{16,8}$
is injective.  Since the $Q_i$'s take integral elements to integral
elements, we have an inclusion
\[ Q_2Q_1 \colon \intH^{6,4}(\CQa;\Z/2) \inc \intH^{16,8}(\CQa;\Z/2).
\]
But it follows directly from Proposition~\ref{pr:rost} that
$\intH^{16,8}(\CQa;\Z/2)=0$, and so we are done. 

The argument for general $n$ follows exactly this pattern: one uses
the composite of the operations $Q_1,Q_2,\ldots,Q_{n-2}$, but
everything else is the same.  

\subsection{Summary}
Here is a list of some of the key elements of the proof:
\begin{enumerate}[(1)]
\item The re-interpretation of the Milnor conjecture as a comparison of
different bi-graded motivic cohomology theories.  An extensive
knowledge about such theories allows
one to deduce statements for any smooth simplicial scheme
from statements only about fields (cf. Proposition~\ref{pr:reduce}).
\item Choice of the splitting variety $\Qa$ (needed for
Propositions~\ref{pr:margolis} and \ref{pr:rost}).
\item The introduction and use of \CCech complexes.
\item The construction of Steenrod operations on motivic cohomology
and development of their basic properties, leading to the proof of
Proposition~\ref{pr:margolis}.
\item The `geometric' results of Rost on motives of quadrics, which lead to
Proposition~\ref{pr:rost}.
\end{enumerate}

\subsection{A notable consequence}
\label{re:Hpt}
The integral motivic cohomology groups of a point $H^{p,q}(\Spec F)$
are largely unknown---the exception is when $q=0,1$.  However, the
proof of the Milnor conjecture tells us exactly what $H^{p,q}(\Spec
F;\Z/2)$ is.  First of all, independently of the  Milnor conjecture it
can be shown to vanish when $p\geq q$ and when $q<0$.  By
Proposition~\ref{pr:reduce} (noting that we now know the hypothesis to
be satisfied for all $n$), it follows that
\[ H^{p,q}(\Spec F;\Z/2) \ra H^{p}_{et}(\Spec F;\mu_2^{\tens q}) \]    
is an isomorphism when $p\leq q$ and $q\geq 0$.  As $\mu_2^{\tens
q}\iso \mu_2$, the \'etale cohomology groups are periodic in $q$; that
is, $H^{*}_{et}(\Spec F;\mu_2^{\tens *})\iso
H^*_{Gal}(F;\Z/2)[\tau,\tau^{-1}]$ where $\tau$ has degree $(0,1)$.

The conclusion is that $H^{*,*}(\Spec F;\Z/2)\iso
H^*_{Gal}(F;\Z/2)[\tau]$, where $\tau$ is the canonical class in
$H^{0,1}$ and the Galois cohomology is regarded as the subalgebra
lying in degrees $(k,k)$.  Of course the Milnor conjecture tells us
that the Galois cohomology is the same as mod $2$ Milnor $K$-theory,
and so we can also write 
$H^{*,*}(\Spec F;\Z/2)\iso \bigl ( K^M_*(F)/2 \bigr )[\tau]$.

\subsection{Further reading}
Both the original papers of Voevodsky \cite{Vmilnor,Vmilnor2} are very
readable, and remain the best sources for the proof.  Summaries have
also been given in \cite{M1} and \cite{Su}.  A proof of the general
Bloch-Kato conjecture was recently given in \cite{Vbkato}---the proof
is similar in broad outline to the 2-primary case we described here,
but with several important differences.  See
the introduction to \cite{Vbkato}.

\enlargethispage{0.2in}

Of course in this section I have completely avoided discussing the two
main elements of the proof, namely Propositions~\ref{pr:margolis} and
\ref{pr:rost}.  The proof of Proposition~\ref{pr:margolis} is in
\cite{Vmilnor,Vmilnor2} and is written in a way that can be understood
by most homotopy theorists.  Proposition~\ref{pr:rost} depends on
results of Rost, which seem to be largely unpublished.  See
\cite{R1,R2} for summaries.

For more about why \CCech complexes arise in the proof, see
Proposition~\ref{pr:cechseq} in the next section.

\vfill\eject


\section{Proof of the conjecture on quadratic forms}

In this section and the next I will discuss two proofs of Milnor's
conjecture on quadratic forms.  The first is from \cite{OVV}, the
second was announced in \cite{M2}.  Both depend on Voevodsky's proof
of the norm residue conjecture.  As I keep saying, I'm only going to
give a vague outline of how the proofs go, but with references for
where to find more information on various aspects.  The present
section deals with the \cite{OVV} proof.
  
\subsection{Preliminaries}
Recall that we are concerned with the map $\nu\colon K^M_*(F)/2 \ra
\Gr_I W(F)$ defined by
$\nu(\brn{a_1,\ldots,a_n})=\br{1,-a_1}\cdots\br{1,-a_n}$.  The fact
that $I$ is additively generated by the forms $\br{1,x}$ shows that
$\nu$ is obviously surjective; so our task is to prove injectivity.
In general, the product $\br{1,b_1}\cdots\br{1,b_n}$ is called an
\dfn{\mdfn{$n$}-fold Pfister form}, and denoted
$\dbr{b_1,\ldots,b_n}$.  Note that it has dimension $2^n$.  The proof
is intimately tied up with the study of such forms.

Milnor proved that the map $\nu\colon K^M_2(F)/2 \ra I^2/I^3$ is an
isomorphism.  He used ideas of Delzant \cite{De} to define
Stiefel-Whitney invariants for quadratic forms, which in dimension $2$
give a map $I^2/I^3 \ra K^M_2(F)/2$.  One could explicitly check that
this was an inverse to $\nu$.  Unfortunately, this last statement
generally fails in larger dimensions; the Stiefel-Whitney invariants
don't carry enough information.  See \cite[4.1, 4.2]{Mi}.

\subsection{The Orlov-Vishik-Voevodsky proof}
We first need to recall some results about Pfister forms proven in the
70's.  The first is an easy corollary of the so-called Main Theorem of
Arason-Pfister (cf. \cite[4.5.6]{S1}).  For a proof, see
\cite[pp. 192-193]{EL}. 
 
\begin{prop}[Elman-Lam]
$\dbr{a_1,\ldots,a_n} \equiv \dbr{b_1,\ldots,b_n} \  (\mmod I^{n+1})$ if
and only if $\dbr{a_1,\ldots,a_n}=\dbr{b_1,\ldots,b_n}$ in $GW(F)$.
\end{prop}

Combining the result for $n=2$ with Milnor's theorem that $K^M_2(F)
\ra I^2/I^3$ is an isomorphism, we get the following
(note that the minus signs are there because
$\nu(\brn{a_1,\ldots,a_n})=\dbr{-a_1,\ldots,-a_n}$):

\begin{cor}
$\dbr{a_1,a_2} = \dbr{b_1,b_2}$ in $GW(F)$
if
and only if $\brn{-a_1,-a_2}=\brn{-b_1,-b_2}$ in $K^M_*(F)/2$.
\end{cor}

Say that two $n$-fold Pfister forms $A=\dbr{a_1,\ldots,a_n}$ and
$B=\dbr{b_1,\ldots,b_n}$ are \dfn{simply-\mdfn{$p$}-equivalent} if
there are two indices $i$, $j$ where $\dbr{a_i,a_j}=\dbr{b_i,b_j}$ and
$a_k=b_k$ for all $k\notin \{i,j\}$.  The forms $A$ and $B$ are
\dfn{chain-\mdfn{$p$}-equivalent} if there is a chain of forms
starting with $A$ and ending with $B$ in which every link of the chain
is a simple-$p$-equivalence.  Note that it follows immediately from
the previous corollary that if $A$ and $B$ are chain-$p$-equivalent
then $\brn{-a_1,\ldots,-a_n}=\brn{-b_1,\ldots,-b_n}$.

The following result is \cite[Main Theorem 3.2]{EL}:

\begin{prop}
Let $A=\dbr{a_1,\ldots,a_n}$ and $B=\dbr{b_1,\ldots,b_n}$.
The following are equivalent:
\begin{enumerate}[(a)]
\item $A$ and $B$ are chain-$p$-equivalent.
\item $\brn{-a_1,\ldots,-a_n}=\brn{-b_1,\ldots,-b_n}$ in $K^M_*(F)/2$.
\item $A\equiv B\ (\mmod I^{n+1})$.
\item $A=B$ in $GW(F)$.
\end{enumerate}
\end{prop}

Note that $(a)\Rightarrow (b) \Rightarrow (c)$ is trivial, and
$(c)\Rightarrow (d)$ was mentioned above.  So the new content is in
$(d)\Rightarrow (a)$.  I will not give the proof, but refer the reader
to \cite[4.1.2]{S1}.  The result below is a
restatement of $(c)\Rightarrow (b)$:

\begin{cor}
\label{co:pureinject}
The equality $\nu(\brn{a_1,\ldots,a_n})=\nu(\brn{b_1,\ldots,b_n})$ 
can only occur if
$\brn{a_1,\ldots,a_n}=\brn{b_1,\ldots,b_n}$.  
\end{cor}

Unfortunately the above corollary does not show that $\nu$ is
injective, as a typical element $x\in K^M_*(F)/2$ is a {\it sum\/} of
terms $\brn{a_1,\ldots,a_n}$.  A term $\brn{a_1,\ldots,a_n}$ is called
a \dfn{pure symbol}, whereas a general $x\in K^M_*(F)$ is just a
\dfn{symbol}.  The key ingredient needed from \cite{OVV} is the
following:

\begin{prop}
\label{pr:expand}
If $x\in K^M_*(F)/2$ is a nonzero element then there is a field
extension $F\inc F'$ such that the image of $x$ in $K^M_*(F')/2$ is a
nonzero pure symbol.
\end{prop}

It is easy to see that the previous two results prove the injectivity
of $\nu$.  If $x\in K^M_n(F)/2$ is a nonzero element in the kernel of
$\nu$, then by passing to $F'$ we find a nonzero pure symbol which is
also in the kernel.  Corollary~\ref{co:pureinject} shows this to be
impossible, however.

We are therefore reduced to proving Proposition~\ref{pr:expand}.  If
we write $x=\und{a}_1+\ldots+\und{a}_k$, where each $\und{a}_i$ is a
pure symbol, then we know we can make $\und{a}_i$ vanish by passing to
the function field $F(Q_{\und{a}_i})$ (where $Q_{\und{a}_i}$ is the
splitting variety produced in the last section).  Our goal will be to
show that $\und{a}_i$ is the {\it only\/} term that vanishes:

\begin{prop}[Orlov-Vishik-Voevodsky]
\label{pr:OVV}
If $\und{a}=\brn{a_1,\ldots,a_n}$ is nonzero in $K^M_n(F)/2$, then the kernel
of $K^M_n(F)/2 \ra K^M_n(F(\Qa))/2$ is precisely $\Z/2$ (generated by
$\und{a}$).
\end{prop}

Granting this for the moment, let $i$ be the largest index for which
$x$ is nonzero in $K^M_n(F')/2$, where
$F'=F(Q_{\und{a}_1}\times\cdots\times Q_{\und{a}_i})$.  Since $x$ will
become zero over $F'(Q_{\und{a}_{i+1}})$, the above result says that
$x=\und{a}_{i+1}$ in $K^M_n(F')/2$.  This is precisely what we wanted.

So finally we have reduced to the same kind of problem we tackled in
the last section, namely controlling the map $K^M_n(F)/2 \ra
K^M_n(F(\Qa))/2$.  The techniques needed to prove
Proposition~\ref{pr:OVV} are exactly the same as those from the last
section.  There is a again a homotopical ingredient and a geometric
ingredient.  

\begin{prop}
\label{pr:cechseq}
If $X$ is a smooth scheme over $F$, then for every $n\geq 0$ there is
an exact sequence of the form
\[ 0 \ra H^{n,n-1}(\Cech X;\Z/2) \ra H^{n,n}(\Spec F;\Z/2)\ra H^{n,n}(\Spec
F(X);\Z/2).
\]
\end{prop}

Recall that $H^{n,n}(\Spec E;\Z/2)\iso K^M_n(E)/2$ for any field $E$.  So
the above sequence is giving us control over the kernel of $K^M_*(F)/2
\ra K^M_*(F(\Qa))/2$.  The proof uses the conclusion from 
Proposition~\ref{pr:reduce} (which is known by Voevodsky's proof of
the Milnor conjecture) and some standard manipulations with motivic
cohomology.  See \cite[Prop. 2.3]{OVV}.

\begin{remark}
\label{re:cech}
In some sense Proposition~\ref{pr:cechseq} explains why \CCech
complexes are destined to come up in the proofs of these conjectures.
\end{remark}

If the above proposition is thought of as a `homotopical' part of
the proof, the geometric part is the following.  It is deduced using
Rost's results on the motive of $\Qa$; see \cite[Prop. 2.5]{OVV}.

\begin{prop}
\label{pr:rost2}
There is a surjection $\Z/2 \ra H^{2^n-1,2^{n-1}-1}(\Cech \Qa;\Z/2)$.
\end{prop}  

The previous two results immediately yield a proof of \ref{pr:OVV}.
By Proposition~\ref{pr:cechseq} we must show that $H^{n,n-1}(\Cech
\Qa;\Z/2)\iso \Z/2$ (and we know the group is nontrivial).  But we saw
in the last section that $H^{n,n-1}(\Cech \Qa;\Z/2)\iso
H^{n+1,n-1}(\tC \Qa;\Z/2)$, where $\tC \Qa$ is the homotopy cofiber of
$(\Cech \Qa)_+ \ra (\Spec F)_+$.  We also saw that the operation
$Q_{n-2}\cdots Q_2Q_1$ gives a monomorphism $H^{n+1,n-1}(\tC \Qa;\Z/2)
\inc H^{2^n-1,2^{n-1}-1}(\tC \Qa;\Z/2)$.  But now we are done, since
by \ref{pr:rost2} the latter group has at most two elements.  

This completes the proof of the injectivity of $\nu$.

\vfill\eject


\section{Quadratic forms and the Adams spectral sequence}
In \cite{M2} Morel announced a proof of the quadratic form conjecture
over characteristic zero fields, using the motivic Adams spectral
sequence.  The approach depends on having computed the motivic
Steenrod algebra, but I'm not sure what the status of this
is---certainly no written account is presently available.  Despite this
frustrating point, Morel's proof is very exciting; while it uses
Voevodsky's computation of $H^{*,*}(\Spec F;\Z/2)$---see
Remark~\ref{re:Hpt}---it somehow avoids using any other deep results
about quadratic forms!  So I'd like to attempt a sketch.

The arguments below take place in the motivic stable homotopy
category.  All the reader needs to know as background is that it
formally behaves much as the usual stable homotopy category, and that
there is a bigraded family of spheres $S^{p,q}$.  The suspension (in
the triangulated category sense) of $S^{p,q}$ is $S^{p+1,q}$, and
$S^{2,1}$ is the suspension spectrum of the variety $\P^1$.

\vspace{0.1in}

\subsection{Outline}
We have our maps $\nu_n\colon K^M_n(F)/2 \ra I^n/I^{n+1}$, and need to
prove that they are injective.  We will see that the Adams spectral
sequence machinery gives us, more or less for free, maps $s_n\colon
I^n/I^{n+1} \ra K^M_n(F)/(2,J)$ where $J$ is a subgroup of boundaries
from the spectral sequence.  The composite $s_n\nu_n$ is the natural
projection, and so the whole game is to show that $J$ is zero.  That
is, one needs to prove the vanishing of a line of differentials.
Using the multiplicative structure of the spectral sequence and the
algebra of the $E_2$-term, this reduces just to proving that the
differentials on a certain `generic' element vanish.  This allows
one to reduce to the case of the prime field $\Q$, then to $\R$, and
ultimately to a purely topological problem.

\vspace{0.1in}

\subsection{Basic setup}
Now I'll expand on this general outline.  The first step is to produce
a map $q\colon GW(F) \ra \{S^{0,0},S^{0,0}\}$ where
$\{\blank,\blank\}$ denotes maps in the motivic stable homotopy
category.  Recall from Section~\ref{se:GW} that one knows a complete
description of $GW(F)$ in terms of generators and relations.  For
$a\in F^*$ we let $q(\br{a})$ be the map $\P^1 \ra \P^1$ defined in
homogeneous coordinates by $[x,y] \ra [x,ay]$.  By writing down
explicit $\A^1$-homotopies one can verify that the relations in
$GW(F)$ are satisfied in $\{S^{0,0},S^{0,0}\}$, and so $q$ extends to
a well-defined map of abelian groups.  It is actually a ring map.
Further details about all this are given in \cite{M3}.

Now we build an Adams tower for $S^{0,0}$ based on the motivic
cohomology spectrum $H\Z/2$.  
Set $W_0=S^{0,0}$, and define $W_1$ by the homotopy fiber sequence
$W_1 \ra S^{0,0}\ra H\Z/2$.  Then consider the map $W_1 \iso
S^{0,0}\Smash W_1 \ra H\Z/2 \Smash W_1$, and let $W_2$ be the homotopy
fiber.  Repeat the process to define $W_3$, $W_4$, etc.  This
gives us a tower of cofibrations 
\[
\xymatrix{
&H \Smash W_2 & H \Smash W_1 & H\Smash W_0 \\
\cdots \ar[r] & W_2\ar[u] \ar[r] & W_1\ar[u] \ar[r]  & W_0,\ar[u]
}
\]
where we have written $H$ for $H\Z/2$.  For any $Y$ the tower yields a
filtration on $\{Y,S^{0,0}\}$ by letting $\cF^n$ be the subgroup of
all elements in the image of $\{Y,W_n\}$ (note that there is no {\it a
priori\/} guarantee that the filtration is Hausdorff.)  The tower
yields a homotopy spectral sequence whose abutment has something to do
with the associated graded of the groups $\{S^{*,0}\Smash
Y,S^{0,0}\}$.  If the filtration is not Hausdorff these associated
graded groups may not be telling us much about $\{S^{*,0}\Smash
Y,S^{0,0}\}$, but this will not matter for our application.  We will
be interested in the case $Y=S^{0,0}$.  

Set $E_1^{a,b}=\{S^{a,0},H\Smash W_b\}$, so that $d_r\colon E_r^{a,b}
\ra E_r^{a-1,b+r}$.  My indexing has been chosen so that the picture
of the spectral sequence has $E_1^{a,b}$ in spot $(a,b)$ on a grid,
rather than at spot $(b-a,a)$ as is more typical for the Adams
spectral sequence---but the picture itself is the same in the end.
Formal considerations give inclusions
\[ \cF^k\{S^{n,0},S^{0,0}\}/\cF^{k+1}\{S^{n,0},S^{0,0}\} \inc E_\infty^{n,k}
\]
(however, there is no a priori reason to believe the map is surjective).
In particular, if $\cF^*$ is the filtration on
$\{S^{0,0},S^{0,0}\}$ then we have inclusions $\cF^k/\cF^{k+1} \inc
E_\infty^{0,k}$.

Let $\Gi(F)$ be the kernel of the mod $2$ dimension function $\dim
\colon GW(F) \ra \Z/2$.  The powers $\Gi(F)^n$ define a filtration on
$GW(F)$.  One can check that $q$ maps $\Gi^1$ into $\cF^1$.
Since the Adams filtration $\cF^n$ on $\pi_{0,0}(S^{0,0})$ will be
multiplicative, one finds that $q$ maps $\Gi^n$ into $\cF^n$.  So
we get maps $(\Gi)^n/(\Gi)^{n+1} \ra \cF^n/\cF^{n+1} \ra E_\infty^{0,n}$.  

In a moment I'll say more about what the Adams spectral sequence looks
like in this case, but first let's relate $\Gi$ to what we really care
about.  One easily checks that $\Gi=GI \oplus \Z$, where the $\Z$ is
the subgroup generated by $\br{1,1}=2\br{1}$.  So $\Gi^n= GI^n \oplus
\Z$, where the $\Z$ is generated by $2^n\br{1}$.  It follows that
$\Gi^n/\Gi^{n+1} \iso [GI^n/GI^{n+1}] \oplus \Z/2$.  Finally,
recall from Remark~\ref{re:GI} that the natural map $GI \ra I$ is an
isomorphism.  Putting everything together, we have produced invariants
$[I^n/I^{n+1}] \oplus \Z/2 \ra E_\infty^{0,n}$.

\vspace{0.1in}

\subsection{Analysis of the spectral sequence}
So far the discussion has been mostly formal.  We have produced a
spectral sequence, but not said anything concrete about it.  The
usefulness of the above invariants hinges on what $E_\infty^{0,n}$
looks like.  If things work as in ordinary topology, then the $E_2$
term will turn out to be
$E_2^{a,b}=\Ext_{H^{**}H}^{b}(\Sigma^{b+a,0}H^{**},H^{**})$ where I've
again written $H=H\Z/2$ and $\Sigma^{k,0}$ denotes a grading shift on
the bi-graded module $H^{**}$. 
So we need to know the algebra $H^{**}H$, but
unfortunately there is no published source for this calculation.  In
\cite{V2} Voevodsky defines Steenrod operations and shows that they
satisfy analogs of the usual Adem relations; he doesn't show that
these generate all of $H^{**}H$, though.  However, let's assume we
knew this---so we are assuming $H^{**}H$ is the algebra Voevodsky
denotes $A^{**}$ and calls the motivic Steenrod algebra \cite[Section
11]{V2}.

The form of $H^{**}H$ is very close to that of the usual Steenrod
algebra, and so one has a chance at doing some of the $\Ext$
computations.  In fact, it is not very hard.  Some hints about this
are given in Appendix B, but for now let me just tell you the
important points:

\begin{enumerate}[(1)]
\item $E_2^{p,q}=0$ if $p<0$.
\item $E_2^{0,0}=\Z/2$.
\item For $n\geq 1$, $E_2^{0,n}=H^{n,n} \oplus \Z/2$.  The inclusion $\oplus_n
H^{n,n} \inc \oplus_n E_2^{0,n}$ is a ring homomorphism, where the
domain is regarded as a subring of $H^{**}$.
\end{enumerate}
Most of these computations make essential use of Remark~\ref{re:Hpt},
and therefore depend on Voevodsky's proof of the norm residue
conjecture.  Also note the connection between (3) and Milnor
$K$-theory, given by the isomorphism $H^{n,n}\iso K^M_n(F)/2$.

The above two facts show that everything in $E_2^{0,n}$ is a permanent
cycle and thus $E_\infty^{0,n}=(\Z/2 \oplus K^M_n(F)/2)/J$ where $J$
is the subgroup of all boundaries.  Recall that one has maps
\[ K^M_n(F)/2 \llra{\nu_n} I^n/I^{n+1} \ra E_\infty^{0,n}\iso 
[K^M_n(F)/2 \oplus \Z/2]/J.
\]
The composition can be checked to be the obvious one.  To prove that
$\nu_n$ is injective, we need to prove that $J=0$.  That is, we need
to prove the vanishing of all differentials landing in $E^{0,*}$
(which necessarily come from $E^{1,*}$).  As for the computation of
the $E^{1,*}$ column, here are the additional facts we need:
\begin{enumerate}[(1)]
\addtocounter{enumi}{3}
\item $E_2^{1,0}=0$.
\item $E_2^{1,1}=H^{0,1}\oplus H^{2,2}\iso \Z/2 \oplus H^{2,2}$.
\item The images of the two maps
\[ E_2^{0,1} \tens E_2^{1,n-1} \ra E_2^{1,n} \qquad E_2^{1,n-1} \tens
E_2^{0,1} \ra E_2^{1,n}
\]
generate $E_2^{1,n}$ as an abelian group.
\item The composite $H^{1,1} \tens H^{2,2} \inc E_2^{0,1} \tens
E_2^{1,1} \ra E_2^{1,2}$ is zero.
\end{enumerate}

Again, let me say that none of these computations is particularly
difficult, and the reader can find some hints in Appendix B.  Portions
of columns $0$ and $1$ of our $E_2$-term are
shown below:

\begin{picture}(340,180)(-30,-10)
\multiput(0,0)(70,0){5}{\line(0,1){150}}
\multiput(0,0)(0,30){6}{\line(1,0){280}}
\put(29,10){$\Z/2$}
\put(8,40){$H^{1,1}\oplus \Z/2$}
\put(8,70){$H^{2,2}\oplus \Z/2$}
\put(8,100){$H^{3,3}\oplus \Z/2$}
\put(8,130){$H^{4,4}\oplus \Z/2$}
\put(99,10){$0$}
\put(78,40){$H^{2,2}\oplus \Z/2$}
\put(99,70){$??$}
\put(99,100){$??$}
\put(99,130){$??$}
%
\put(-2,-2){$\bullet$}
\thicklines
\put(0,0){\vector(0,1){155}}
\put(0,0){\vector(1,0){285}}
\end{picture}

\begin{remark}
If one only looks at the $\Z/2$'s appearing in the above diagram, the
picture looks just like the ordinary topological Adams spectral
sequence.  The $\Z/2$'s in our $0$th column indeed turn out to be
``$h_0^n$'s'', just as in topology.  The $\Z/2$ in $E_2^{1,1}$ is a
little more complicated, though---it doesn't just come from $\Sq^2$,
like the usual $h_1$ does (see Appendix B for what it {\it does\/}
come from).
\end{remark}

We need to prove that all the differentials leaving the $E^{1,*}$
column vanish.  By fact (6) and the multiplicative structure of the
spectral sequence, it is sufficient to prove that all differentials
leaving $E_2^{1,1}$ vanish (starting with $d_2\colon
E_2^{1,1} \ra E_2^{0,3}$).  We will do this in several steps.

The following result basically shows that, just as in ordinary
topology, all the $\Z/2$'s in column $0$ survive to $E_\infty$.

\begin{lemma}
The image of $d_r \colon E_r^{1,1} \ra E_r^{0,r+1}$
lies in the subgroup $H^{r+1,r+1}$, for every $r\geq 2$.
\end{lemma}

\begin{proof}
Suppose there is an element $x\in E_r^{1,1}$ such that $d_r(x)$ does
not lie in $H^{r+1,r+1}$ (or rather its image in $E_r$).  We can write
$x= \und{a}+y$ where $\und{a} \in H^{2,2}=K^M_2(F)/2$ and $y\in
H^{0,1}\iso \Z/2$.  In expressing $\und{a}$ as a sum of pure symbols,
one notes that only a finite number of elements of $F$ are involved.
By naturality of the spectral sequence, we can therefore assume $F$ is
a finitely-generated extension of $\Q$.

But now we can choose an embedding $F\inc \C$, and again use
naturality.  The groups $K^M_n(\C)/2$ are all zero, and therefore our
assumption implies that over $\C$ we have $E_{r+1}^{0,r+1}=0$ (in
other words, the $\Z/2$ in $E_2^{0,r+1}$ dies in the spectral
sequence).  But there is a `topological realization map' from our
spectral sequence over $\C$ to the usual Adams spectral sequence in
topology, where we know that none of the $\Z/2$'s in $E^{0,*}$ ever die.
\end{proof}

\begin{remark}
\label{re:algprf}
There is also a purely algebraic proof of the above result.  One
reduces via naturality to the case of algebraically closed fields,
where all the $H^{n,n}$'s are zero.  Then one shows that the $\Z/2$'s
in the $0$th column form a polynomial algebra, and that the composite
$\Z/2 \tens \Z/2 \inc E_2^{1,1} \tens E_2^{0,1} \ra E_2^{1,2}$ is zero
(just as in ordinary topology).  The fact that the spectral sequences
is multiplicative takes care of the rest.
\end{remark}

\begin{lemma}
For $\und{a}\in H^{2,2}$ one has $d_r(\und{a})=0$, for every $r$.
\end{lemma}

\begin{proof}
It follows from facts (3) and (7), together with the multiplicative
structure of the spectral sequence, that everything in the image of
$d_r \colon H^{2,2} \ra H^{r+1,r+1}$ is killed by $H^{1,1}$.  This is
the key to the proof.

Let $z=d_r(\und{a})$.  Consider the naturality of the spectral
sequence for the map $j\colon F \ra F(t)$.  It follows from the
previous paragraph that $j(z)=d_r(j\und{a})$ is killed by $F(t)^*$.  
In particular, $\brn{t} \cdot j(z) =0$ in $K^M_{r+2}(F(t))/2$.  But by
\cite[Lem. 2.1]{Mi} there is a map $\partial_t \colon K^M_{r+2}(F(t))/2
\ra K^M_{r+1}(F)/2$ with the property that $\partial_t(\brn{t}\cdot
j(z))=z$.  So we conclude that $z=0$, as desired. 
\end{proof}

\begin{prop}
All differentials leaving $E^{1,1}$ are zero.
\end{prop}

\begin{proof}
Recall $E_2^{1,1}\iso H^{0,1}\oplus H^{2,2}\iso \Z/2 \oplus H^{2,2}$.
By the previous lemma we are reduced to analyzing the maps $d_r\colon
H^{0,1} \ra H^{r+1,r+1}$.  Since $H^{0,1}(\Q) \ra H^{0,1}(F)$ is an
isomorphism, it suffices to prove the result in the case $F=\Q$.

Now use naturality with respect to the field extension $\Q \inc \R$.
The maps $K^M_n(\Q)/2 \ra K^M_n(\R)/2$ are isomorphisms for $n\geq 3$
(see Appendix A), so now we've reduced to $F=\R$.  But here we can
again use a `topological realization' map to compare our Adams
spectral sequence to the corresponding one in the context of
$\Z/2$-equivariant homotopy theory.  This map is readily seen to be an
isomorphism on the $E^{0,*}$ column: the point is that the
$\Z/2$-equivariant cohomology groups $H^{n,n}$ are isomorphic to the
corresponding mod $2$ motivic cohomology groups over $\R$ (see
\cite[2.8, 2.11]{Du}, for instance).  We are essentially seeing a
reflection of the fact that $GW(\R)$ may be identified with the
Burnside ring of $\Z/2$, which coincides with $\{S^{0,0},S^{0,0}\}$ in
the $\Z/2$-equivariant stable homotopy category.  In any case, we are
finally reduced to showing the vanishing of certain differentials in a
topological Adams spectral sequence: the paper \cite{LZ} seems to
essentially do this (but I haven't thought about this part
carefully---I'm relying on remarks from \cite{M2}).
\end{proof}

This completes Morel's proof of the quadratic form conjecture for
characteristic zero fields (modulo the identification of $H^{**}H$,
which we assumed).

\vspace{0.2in}

\begin{remark}
We restricted to characteristic zero fields because the identification
of $H^{**}H$ has never been claimed in characteristic $p$.  If we make
the wild guess that in positive characteristic $H^{**}H$ still has the
same form, most of the argument goes through verbatim.  There are two
exceptions, where we used topological realization functors.  The first
place was to show that the image of the $d_r$'s didn't touch the
$\Z/2$'s in $E_2^{0,*}$, but Remark~\ref{re:algprf} mentioned that
this could be done another way.  The second place we used
topological realization was at the final stage of the argument, to
analyze the differentials $d_r \colon H^{0,1} \ra H^{r+1,r+1}$.  As
before, this reduces to the case of a prime field.  But for $F$ a
finite field one has $K^M_n(F)=0$ for $n\geq 2$, so for prime fields
there is in fact nothing to check.

In summary, the same general argument would work in characteristic
$p$ if one knew that $H^{**}H$ had the same form.
\end{remark}

\subsection{Further reading}
There is very little completed literature on the subjects discussed in
this section.  Several documents are available on Morel's website,
however; the draft \cite{M5} is particularly relevant, although it
only slightly expands on \cite{M2}.  For information on the motivic
Steenrod algebra, see \cite{V2}.  Finally, Morel recently released
another proof of Milnor's quadratic form conjecture, using very
different methods.  See \cite{M4}.  

\vfill\eject


\appendix

\section{Some examples of the Milnor conjectures}

This is a supplement to Section 1.  We examine the Milnor conjectures
in the cases of certain special fields $F$.

\vspace{0.1in}

\subsubsection{(a) $F$ is algebraically closed.}
Since $F=F^2$, every nondegenerate form is isomorphic to one of the
form $\langle 1,1,\ldots,1 \rangle$.  So $GW(F)\iso \Z$, and $W(F)\iso
\Z/2$ with $I(F)=0$.  Thus, $\Gr_I W(F)\iso \Z/2$.

The absolute Galois group is trivial, so $H^*(F;\Z/2)=\Z/2$.

Finally, the fact that $F=F^2$ implies that $K^M_*(F)/2=0$ for $*\geq
1$.  This is because the generators all lie in $K^M_1(F)$, and if
$a=x^2$ then $\brn{a}=\brn{x^2}=2\brn{x}=0\in K^M_1(F)/2$.

\vspace{0.1in}

\subsubsection{(b) $F=F^2$.}
This case is suggested by the previous one.  We only need to check
that the hypothesis implies $H^*(F;\Z/2)=0$ for $*\geq 1$.  Strangely,
I haven't been able to find an easy proof of this.

\vspace{0.1in}

\subsubsection{(c) $F=\R$.}

In this case we know forms are classified by their rank and signature,
and it follows that $GW(\R)$ is the free abelian group generated by
$\br{1}$ and $\br{-1}$. Also, $\br{-1}^2=\br{1}$.  So $GW(\R)\iso
\Z[x]/(x^2-1)$, and $W(\R) \iso \Z$ with $I(\R)=2\Z$.  Hence $\Gr_I
W(\R)\iso \Z/2[a]$.

The absolute Galois group of $\R$ is $\Z/2$, so
$H^*(\R;\Z/2)=H^*(\Z/2;\Z/2)=\Z/2[a]$.  

Finally we consider $K^M_*(\R)/2$.  The group
$K^M_1(\R)/2=\R^*/(\R^*)^2 \iso \{1,-1\}$ (the set consisting of $1$
and $-1$).  A similar calculation, based on the fact that every
element of $\R$ is a square up to sign, shows that $K^M_i(\R)/2\iso
\Z/2$ for every $i$, with the nonzero element being
$\brn{-1,-1,\ldots,-1}$.  So $K^M_*(\R)/2 \iso \Z/2[a]$ as well.

\vspace{0.1in}

\subsubsection{(d) $F=\F_q$, $q$ odd.}
Here $F^*\iso \Z/(q-1)$ and so $K^M_1/2=F^*/(F^*)^2 \iso \Z/2$.  If
$g$ is the generator, then $\brn{g,g,\ldots,g}$ generates $K^M_n/2$
(but may be zero).  In fact one can show (cf. \cite[Ex. 1.5]{Mi}) that
$\brn{g,g}=0$ in $K^M_2$, from which it follows that $K^M_*=0$ for
$*\geq 2$.  So $K^M_*(F)/2\iso \Z/2\oplus \Z/2$, in degrees $0$ and
$1$.

For a finite field the absolute Galois group is $\hat{\Z}$, the
profinite completion of $\Z$.  The Galois cohomology
$H^*(\hat{\Z};\Z/2)$ is just the mod $2$ cohomology of $B\Z \he S^1$;
so it is $\Z/2\oplus \Z/2$, with the generators in degrees $0$ and
$1$.

Again noting that $F^*/(F^*)^2\iso \Z/2$, it follows that 
the Grothendieck-Witt group is generated by $\br{1}$ and $\br{g}$.
A simple counting argument (cf. \cite[Lem. 2.3.7]{S1}) shows that
every element of $\F_q^*$ is a sum of two squares.  Writing $g=a^2+b^2$
one finds that 
\[
\br{1,1}=\br{a^2,b^2}=\br{a^2+b^2,a^2b^2(a^2+b^2)}=\br{a^2+b^2,a^2+b^2}
=\br{g,g}.
\]
That is, $2(\br{1}-\br{g})=0$.
It follows that $GW(F)=\Z\oplus \Z/2$, with corresponding generators
$\br{1}$ and $\br{1}-\br{g}$.

The computation of the Witt group depends on whether or not $-1$ is a
square; since $F^*=\Z/(q-1)$ and $-1$ has order $2$, then $-1$ is a
square precisely when $4|(q-1)$.  So if $q\equiv 1 (\mmod 4)$ then
$\br{1}=\br{-1}$ and $W(F)\iso \Z/2\oplus \Z/2$; in this case
$I(F)=(\br{1}-\br{g})\iso \Z/2$.  If $q\equiv 3 (\mmod 4)$ then
$\br{g}=\br{-1}$ and we have $W(F)\iso \Z/4$ with $I(F)=(2)$.  In
either case $\Gr_I W(F)\iso \Z/2\oplus \Z/2$.

\begin{remark}
Although Milnor's quadratic form conjecture says that $\Gr_I W(F)$
depends only on the absolute Galois group of $F$, this example makes
it clear that the same cannot be said for $W(F)$ itself.
\end{remark}

\subsubsection{(e) $F=\Q$.}
This case is considerably harder, so we will only make a few observations.
Note that as an abelian group one has
\[ \Q^* \iso \Z/2 \times \Bigl (\oplus_p \Z \Bigr ),
\]
by the fundamental theorem of arithmetic; the direct sum is over the
set of all primes.  Here the isomorphism sends a fraction $q$ to its
sign (in the $\Z/2$ factor) together with the list of exponents in the
prime factorization of $q$.  So $K^M_1(\Q)/2 \iso \Z/2 \oplus
(\oplus_p \Z/2)$.  

As the above isomorphism may suggest,
to go further it becomes convenient to work with one completion at a
time.  The case $F=\R$ has already been discussed, so what is left is
the $p$-adics.  We will return to $F=\Q$ after discussing them.

\vspace{0.1in}

\subsubsection{(f) $F=\Q_p$.}
We will concentrate on the case where $p$ is odd; the case $p=2$ is
similar, and can be left to the reader.
We know $K^M_1(\Q_p)/2 \iso H^1(\Q_p;\Z/2)\iso \Q_p^*/(\Q_p^*)^2$.  A
little thought (cf. \cite[5.6.2]{S1}) shows this group is $\Z/2\oplus
\Z/2$, with elements represented by $1$, $g$, $p$, and $pg$, where $1<
g < p$ is any integer which generates the multiplicative group
$\F_p^*$.  By \cite[Section II.5.2]{Se} one has $H^2(\Q_p;\Z/2)\iso \Z/2$
and $H^i(\Q_p;\Z/2)=0$ for $i\geq 3$.

The fact that $K_1^M(\Q_p)/2$ only has four elements tells us that
$K^M_*(\Q_p)/2$ can't be too big.  By finding the appropriate
relations to write down, Calvin Moore proved that $K^M_*(\Q_p)/2=0$
for $*\geq 3$ \cite[Ex. 1.7]{Mi}, and that $K^M_2(\Q_p)/2=\Z/2$.
This is an exercise for the reader.

The group $GW(\Q_p)$ will be generated by the four elements $\br{1}$,
$\br{g}$, $\br{p}$, and $\br{pg}$.  The theory again depends on
whether or not $-1$ is a square, which is when $p\equiv 1(\mmod 4)$.  
When $p\equiv 1 (\mmod 4)$ one has $\br{1}=\br{-1}$ and so
$\br{x}=\br{-x}$ for any $x$.  As a result
$\br{g,g}=\br{g,-g}=\br{1,-1}=\br{1,1}$, and similarly
$\br{p,p}=\br{pg,pg}=\br{1,1}$.  One finds that $GW(\Q_p)=\Z\oplus
(\Z/2)^3$ with corresponding generators $\br{1}$, $\br{1}-\br{p}$,
$\br{1}-\br{g}$, and $\br{1}-\br{pg}$.  
Since $\br{1,-1}=2\br{1}$, $W(\Q_p)=(\Z/2)^4$ with the same
generators.  $I$ is generated by $\br{1,p}$, $\br{1,g}$, and
$\br{1,pg}$; $I^2$ is generated by $\br{1,p,g,pg}$, and $I^3=0$.
So $\Gr_I W=\Z/2 \oplus (\Z/2\oplus \Z/2) \oplus \Z/2$.  Note that
this is the first example we've seen where $I^2 \neq 2I$.

When $p\equiv 3 (\mmod 4)$ we can take $g=-1$.  One has
$\br{1,1}=\br{-1,-1}$ by the same reasoning as for $\F_p$ ($-1$ is the
sum of two squares), and so $\br{p,p}=\br{-p,-p}$.  
Note that
\[ 
\br{p,p,p,p}=\br{p,-p,-p,p}=\br{1,-1,-1,1}=\br{1,1,1,1}
\]
and so $4(\br{1}-\br{p})=0$.  Also, 
\[ \br{p,p,p}=\br{p,-p,-p}=\br{1,-1,-p} \quad\text{and}\quad
\br{1,1,1}=\br{1,-1,-1}.
\]
So $3(\br{1}-\br{p})=\br{-1}-\br{-p}$.  Of course $GW(\Q_p)$ is
generated by $\br{1}$, $\br{1}-\br{-1}$, $\br{1}-\br{p}$, and
$\br{1}-\br{-p}$, and the previous computation shows the last
generator is not needed.  So we have a surjective map $\Z \oplus \Z/2
\oplus \Z/4 \ra GW(\Q_p)$ sending the standard generators to
$\br{1}$, $\br{1}-\br{-1}$, and $\br{1}-\br{p}$.  This is readily
checked to be injective once one knows that $\br{1,1} \not\iso
\br{p,p}$.  If these forms were isomorphic it would follow by
reduction mod some power of $p$ that
$\br{1,1}$ was isotropic over some $\F_{p^e}$; that is, we would have
$\br{1,1} \iso \br{1,-1}$.  But we've already computed $GW(\F_{p^e})$,
and know this is not the case.

The Witt ring is $W(\Q_p)\iso \Z/4\oplus \Z/4$ with generators
$\br{1}$ and $\br{1}-\br{p}$.  The ideal $I$ is generated by $2\br{1}$
and $\br{1}-\br{p}$; $I^2$ is generated by $2(\br{1}-\br{p})$;
$I^3=0$.  Again we have $\Gr_I W\iso \Z/2 \oplus (\Z/2\oplus \Z/2)
\oplus \Z/2$.

\vspace{0.1in}

\subsubsection{(g) Return to $F=\Q$.}
Our understanding of the higher Milnor $K$-groups of $\Q$ is based on
passing to the various completions $\Q_p$ and $\R$.  
A computation of Bass and Tate \cite[Lem. A.1]{Mi} gives an exact sequence
\[ 0\ra K^M_2(\Q)/2 \ra K_2^M(\R)/2 \oplus \Bigl( \oplus_p
K_2^M(\Q_p)/2 \Bigr ) \ra \Z/2 \ra 0,
\] 
and we already know $K^M_2(\Q_p)/2 \iso K^M_2(\R)/2 \iso \Z/2$.
A computation of Tate
\cite[Th. A.2, Ex. 1.8]{Mi} shows that for $*\geq 3$ one has
\[ K^M_*(\Q)/2 \iso \oplus_p K^M_*(\Q_p)/2 \oplus K^M_*(\R)/2 \iso
0\oplus \Z/2.
\]

To compute $H^*(\Q;\Z/2)$ we again work one completion at a time.
A theorem of Tate
\cite[Section II.6.3, Th. B]{Se} says that for $i\geq 3$ one has
\[ H^i(\Q;\Z/2) \iso H^i(\R;\Z/2) \times \prod_p H^i(\Q_p;\Z/2) 
\iso H^i(\R;\Z/2) \iso \Z/2.
\]
Our computation of $\Q^*/(\Q^*)^2\iso H^1(\Q;\Z/2)$ shows that the map
$H^1(\Q;\Z/2) \ra H^1(\R;\Z/2) \times \prod_p H^1(\Q_p;\Z/2)$ is
injective.  More of Tate's work \cite[Sec. II.6.3, Th. A]{Se}
identifies the dual of the kernel with the kernel of $H^2(\Q;\Z/2) \ra
H^2(\R;\Z/2) \times (\oplus_p H^2(\Q_p;\Z/2))$---thus, this latter map
is also injective.  Using this, \cite[Sec. II.6.3, Th. C]{Se} gives a short
exact sequence
\[0 \ra H^2(\Q;\Z/2) \ra H^2(\R;\Z/2) \oplus (\oplus_p H^2(\Q_p;\Z/2))
\ra \Z/2 \ra 0.
\]
As we have already remarked that $H^2(\Q_p;\Z/2)=H^2(\R;\Z/2)=\Z/2$, this
completes the calculation of $H^*(\Q;\Z/2)$.

The method for computing the Witt group $W(\Q)$ proceeds similarly by
working one prime at a time.  See \cite[Section 5.3]{S1}.  One has an
isomorphism of groups $W(\Q)\iso \Z \oplus (\oplus_p W(\F_p))$
\cite[Thm. 5.3.4]{S1}.  With enough trouble one can compute $\Gr_I
W(\Q)$, but we will leave this for the reader to consider.

\begin{remark}
Note that the verification of the Milnor conjectures for $F=\Q$ tells
us exactly how to classify quadratic forms over $\Q$ by invariants.
First one needs the invariants over $\R$ (which are just rank and
signature), and then one needs the invariants over each $\Q_p$---but
for $\Q_p$ one has $I^3=0$, and so $p$-adic forms are classified by
the three classical invariants $e_0$, $e_1$, and $e_2$.  These
observations are essentially the content of the classical
Hasse-Minkowski theorem.
\end{remark}

The method we've used above, of working one completion at a time,
works for all global fields; this is due to Tate for Galois
cohomology, and Bass and Tate for $K^M_*$.  In this way one verifies
the Milnor conjecture for this class of fields \cite[Lemma 6.2]{Mi}.
Note in particular that the class includes all finite extensions of
$\Q$.

\vfill\eject

\section{More on the motivic Adams spectral sequence}

This final section is a supplement to Section 4.  I will give some
hints on computing the $E_2$-term of the motivic Adams spectral
sequence, for the reader who would like to try this at home.  The
computations are not hard, but there are several small issues that are
worth mentioning.

\vspace{0.1in}

\subsection{Setting things up}
\label{se:back}
$H^{**}H$ is the algebra of operations on mod $2$ motivic cohomology.
We will write this as $\mA$ from now on.  There is the Bockstein
$\beta \in \mA^{1,0}$ and there are squaring operations $\Sq^{2i} \in
\mA^{2i,i}$.  We set $\Sq^{2i+1}=\beta \Sq^{2i}\in \mA^{2i+1,i}$.
Finally, there is an inclusion of rings $H^{**} \ra \mA$ sending an
element $t$ to the operation left-multiplication-by-$t$.  Under our
standing assumptions about $\mA$ (see Section 4), it is free as a left
$H^{**}$-module with a basis consisting of the admissible
sequences $\Sq^{i_1}\Sq^{i_2}\cdots \Sq^{i_k}$.

There are two main differences between what happens next and what
happens in ordinary topology.  These are:

\medskip

\begin{enumerate}[(a)]
\item The vector space $H^{**}=H^{**}(\pt)$, regarded as a left
$\mA$-module, is nontrivial.
\item The image of $H^{**} \inc \mA$ is not central.
\end{enumerate}

\medskip

The above two facts are connected.  Let $t\in H^{**}$ and let $\Sq$ denote
some Steenrod operation.  It is not true in general that $\Sq(t\cdot
x) = t\cdot \Sq(x)$---instead there is a Cartan formula for the
left-hand side \cite[9.7]{V2}, which involves Steenrod operations on $t$.
So the operations $\Sq \cdot t$ and $t\cdot \Sq$
are not the same element of $\mA$.  There is one notable exception,
which is when all the Steenrod squares vanish on $t$.  This happens
for elements in $H^{n,n}$, for dimension reasons.  So we have

\medskip

\begin{enumerate}[(a)]
\addtocounter{enumi}{2}
\item Every element of $H^{n,n}$ is central in $\mA$.  
\end{enumerate}

\medskip

It is important that we can completely understand $H^{**}$ as an
$\cA$-module.  This will follow from (1) the fact that $H^{**}\iso
\bigl (\oplus_n H^{n,n}\bigr )[\tau]$ (see Remark~\ref{re:Hpt}); (2)
all Steenrod operations vanish on $H^{n,n}$ for dimension reasons; (3)
all $\Sq^i$'s vanish on $\tau$ except for $\Sq^1$, and
$\Sq^1(\tau)=\rho=\brn{-1}\in H^{1,1}$; (4) the Cartan formula.  In
particular we note the following two facts about $H^{**}$, which are
all that will be needed later (the second fact only needs
Remark~\ref{re:Hpt}):

\medskip

\begin{enumerate}[(a)]
\addtocounter{enumi}{3}
\item The map $\Sq^2 \colon H^{n-1,n} \ra H^{n+1,n+1}$ is zero for all
$n\geq 1$.
\item The map $H^{p,q} \tens H^{i,j} \ra H^{p+i,q+j}$ is surjective
for $q \geq p \geq 0$ and $j\geq i \geq 0$.
\end{enumerate}

\medskip

We are aiming to compute $\Ext^a_{\cA}(H^{**},\Sigma^{b,0}H^{**})$.
In ordinary topology we could use the normalized bar construction to
do this, but one has to be careful here because $H^{**}$, as a left
$\cA$-module, is not the quotient of $\cA$ by a two-sided ideal.  One
way to see this is to use the fact that $\Sq^1(\tau)=\rho$.  Under the
quotient map $\cA \ra H^{**}$ sending $\theta$ to $\theta(1)$, $\Sq^1$
maps to zero but $\Sq^1\tau$ does not (it maps to $\rho$).

So instead of the normalized bar construction we must use the
unnormalized one.  This can be extremely annoying, but for the most
part it turns out not to influence the ``low-dimensional''
calculations we're aiming for.  It is almost certainly an issue when
computing past column two of the Adams $E_2$ term, though.  Anyway,
let
\[ B_n= \mA \tens_{H^{**}} \mA \tens_{H^{**}} \cdots
\tens_{H^{**}} \mA \tens_{H^{**}} H^{**}
\] 
($n+1$ copies of $\mA$).
The final $H^{**}$ 
can be dropped off, of course, but it's useful to keep it there
because the $\mA$-module structure on $H^{**}$ is nontrivial and
enters into the definition of the boundary map.  If we denote the
generators of
$B_n$ as $x=a[\theta_1| \theta_2 | \cdots | \theta_n]t$ then the
differential is
\[ d(x)=(a\theta_1)[\theta_2 | \cdots | \theta_n]t + a[\theta_1\theta_2
| \theta_3 | \cdots | \theta_n]t + \cdots +
  a[\theta_1|\cdots|\theta_{n-1}]\theta_n(t).
\]
The good news is that our coefficients have characteristic $2$, and so
we don't have to worry about signs.  Note that $B_n$, as a left
$H^{**}$-module, is free on generators $1[\theta_1| \cdots
|\theta_n]1$ where each $\theta_i$ is an admissible sequence of
Steenrod operations (and we must include the possibility of the null
sequence $\Sq^0=1$).  We will often drop the $1$'s off of either end
of the bar element, for convenience.

Generators of $\Hom_{\mA}(B_n,H^{**})$ can
be specified by giving a bar element $[\theta_1| \cdots |\theta_n]$
together with an element $t\in H^{**}$.  This data defines a
homomorphism $B_n \ra H^{**}$ sending the generator
$[\theta_1|\cdots|\theta_n]$ to $t$ and all other generators of $B_n$
to zero.  Let's denote this homomorphism by
$t[\theta_1|\cdots|\theta_n]^*$.  These elements generate
$\Hom_{\mA}(B_n,H^{**})$ as an abelian group.  

The last general point to make concerns the multiplicative structure
in the cobar construction.  If we were working with $\Ext_A(k,k)$
where $k$ is commutative and $A$ is an augmented $k$-algebra,
multiplying two of the above generators in the cobar complex just
amounts to concatenating the bar elements---the labels $t\in k$
commute with the $\theta$'s, and so can be grouped together:
e.g. $t\bare{\theta}{n} \cdot u \bare{\alpha}{k} = tu
[\theta_1|\cdots|\theta_n | \alpha_1 |\cdots | \alpha_k]$.  In our
case, the fact that $H^{**}$ is not central in $\mA$ immensely
complicates the product on the cobar complex: very roughly, the $u$
has to be commuted across each $\theta_i$, and in each case a
resulting Cartan formula will introduce new terms into the product.
Luckily there is one case where these complications aren't there,
which is when $u\in H^{n,n}$---for then $u$ is in the center of $\mA$,
and the product works just as above.  We record this observation for
future use:

\medskip

\begin{enumerate}[(a)]
\addtocounter{enumi}{5}
\item $t\bare{\theta}{n}^* \cdot u \bare{\alpha}{k}^*
= tu [\theta_1|\cdots|\theta_n | \alpha_1 |\cdots | \alpha_k]^*$ when
$u\in H^{q,q}$.
\end{enumerate}

\vspace{0.1in}

\subsection{Computations}
We are trying to compute the groups
$\Ext_{\cA}^a(H^{**},\Sigma^{b,0}H^{**})$, and from here on everything
is fairly straightforward.  As an example let's look at $b=1$.  Since
$H^{p,q}\neq 0$ only when $0\leq p\leq q$, one sees that
$\Hom_{\cA}(B_0,H^{**})=0$ and $\Hom_{\cA}(B_1,\Sigma^{1,0}H^{**})\iso
H^{0,0}\oplus H^{1,1}$.  The generators
for this group are elements of the form $s[\Sq^1]^*$ and $t[\Sq^2]^*$,
where $s\in H^{0,0}$ and $t\in H^{1,1}$.  

We likewise find that $\Hom_{\cA}(B_2,\Sigma^{1,0}H^{**})\iso H^{0,1}
\oplus H^{0,1} \oplus H^{0,1} \oplus H^{0,1}$, generated by elements
$s[\Sq^1|1]^*$, $s[1|\Sq^1]^*$, $t[\Sq^2|1]^*$, and $t[1|\Sq^2]^*$.  A
similar analysis shows that $\Hom_{\cA}(B_n,\Sigma^{1,0}H^{**})$ only
has such `degenerate' terms for $n\geq 2$.  No degenerate terms like
these contribute elements to $\Ext$ (at worst they can contribute
relations to $\Ext$).  So the $\Ext^n$'s vanish for $n\geq 2$.  An
analysis of the coboundary shows that everything in dimension $1$ is a
cycle.  So we find that
\[
0=\Ext^0(H^{**},\Sigma^{1,0}H^{**})=\Ext^n(H^{**},\Sigma^{1,0}H^{**}),
\quad \text{for $n\geq 2$}
\]
and
\[ \Ext^1(H^{**},\Sigma^{1,0}H^{**})\iso H^{0,0} \oplus H^{1,1} 
\]
with a typical element in the latter group having the form
$s[\Sq^1]^*+ t[\Sq^2]^*$ (where $s \in H^{0,0}$ and $t\in H^{1,1}$).

\vspace{0.1in}

In general, one sees for degree reasons that the `non-degenerate'
terms in $\Hom_{\cA}(B_n,\Sigma^{n,0}H^{**})$ all have the form
$t\bare{\theta}{n}^*$ where each $\theta_i$ is either $\Sq^1$ or
$\Sq^2$.  In $\Hom_{\cA}(B_{n-1},\Sigma^{n,0}H^{**})$ one has
non-degenerate terms
$u\bare{\theta}{n-1}^*$ of the following types:

\smallskip

\begin{enumerate}[(i)]
\item Each $\theta_i \in \{\Sq^1,\Sq^2\}$, and at least one $\Sq^2$ occurs.
Here $u\in H^{j-1,j}$ where $j$ is the number of $\Sq^2$'s.
\item Each $\theta_i \in \{\Sq^1,\Sq^2,\Sq^3\}$, and exactly one $\Sq^3$
occurs. Here $u\in H^{j+1,j+1}$ where $j$ is the number of $\Sq^2$'s.
\item Each $\theta_i \in \{\Sq^1,\Sq^2,\Sq^2\Sq^1\}$, and exactly one $\Sq^2\Sq^1$
occurs. Here one has $u\in H^{j+1,j+1}$ where $j$ is the number of $\Sq^2$'s.
\item Each $\theta_i \in \{\Sq^1,\Sq^2,\Sq^4\}$, and exactly one $\Sq^4$
occurs. Here $u\in H^{j+2,j+2}$ where $j$ is the number of $\Sq^2$'s. 
\end{enumerate}

\vspace{0.1in}

To analyze the part of the boundary $B_n \ra B_{n-1}$ that we care
about, one only needs to know the Adem relations $\Sq^1\Sq^2=\Sq^3$
and $\Sq^2\Sq^2=\tau \Sq^3\Sq^1$.  (In fact, since $\Sq^3\Sq^1$
doesn't appear in any of the bar elements relevant to
$\Hom(B_{n-1},\Sigma^{n,0}H^{**})$, one may as well pretend
$\Sq^2\Sq^2=0$.)  From this it's easy to compute that
$\Ext^n(H^{**},\Sigma^{n,0}H^{**})\iso H^{0,0}\oplus H^{n,n}$ where a
typical element has the form
$s[\Sq^1|\Sq^1|\cdots|\Sq^1]^*+t[\Sq^2|\Sq^2|\cdots|\Sq^2]^*$.  The
computation uses remark~\ref{se:back}(d).  Also, one sees that all
elements $s[\Sq^1|\Sq^2]^*$ and $s[\Sq^2|\Sq^1]^*$ are zero in
$\Ext^2$ (being the coboundaries of $s[\Sq^3]^*$ and
$s[\Sq^2\Sq^1]^*$, respectively).  Using remark (f) from
Section~\ref{se:back}, this completely determines $\oplus_n
\Ext^n(H^{**},\Sigma^{n}H^{**})$ as a subring of the whole
$\Ext$-algebra.

\vspace{0.2in}

The next step is to compute $\Ext^{0}(H^{**},\Sigma^{1,0}H^{*,*})$,
$\Ext^{1}(H^{**},\Sigma^{2,0}H^{*,*})$, and
$\Ext^{2}(H^{**},\Sigma^{3,0}H^{*,*})$ completely.  The first group is
readily seen to vanish.  For the second group one has to grind out
another term of the bar construction, but it's a very small term.
One finds that
\[
\Ext^{1}(H^{**},\Sigma^{2,0}H^{*,*})\iso H^{0,1} \oplus H^{2,2}
\]
where the generators have the form $s[\Sq^2]^*+(\Sq^1 s)[\Sq^3]^*$ 
and $t[\Sq^4]^*$.  To get the $\Ext^2$ group one will need three more
Adem relations, namely
\[
\Sq^2\Sq^3=\Sq^5+\Sq^4\Sq^1, \quad \Sq^2\Sq^4=\Sq^6+\tau\Sq^5\Sq^1,
\quad \text{and} \quad
\Sq^3\Sq^2=\rho \Sq^3\Sq^1.  
\]
Then the same kind of coboundary
calculations (but a few
more of them) show that
\[ \Ext^{2}(H^{**},\Sigma^{3,0}H^{*,*})\iso H^{1,2} \oplus H^{2,2}
\]
where the generators are $s[\Sq^2|\Sq^2]^*+(\Sq^1 s)[\Sq^3|\Sq^2]^*$
and $t[\Sq^1|\Sq^4]^*=t[\Sq^4|\Sq^1]^*$ (these last two classes are the same
in $\Ext$).  It is important to note that
all elements $u[\Sq^2|\Sq^4]^*$ and $u[\Sq^4|\Sq^2]^*$ are coboundaries
(of $u[\Sq^6]^*$ and $u[\Sq^4\Sq^2]^*$, respectively).  This
justifies fact (7) on page 20.  To justify fact (6) from that same
page (for $n=2$), one notices that the cycles $s[\Sq^2|\Sq^2]^*+(\Sq^1
s)[\Sq^3|\Sq^2]^*$ and $t[\Sq^4|\Sq^1]^*$
decompose as a products
\[ \bigl ( s_1[\Sq^2]^*+(\Sq^1 s_1)[\Sq^3]^* \bigr ) \cdot
(s_2[\Sq^2]^*) \qquad\text{and}\qquad
\bigl ( t_1[\Sq^4]^* \bigr ) \cdot \bigl (t_2[\Sq^1]^* \bigr )
\]
for some $s_1\in H^{0,1}$, $s_2\in H^{1,1}$, $t_1\in H^{2,2}$, and
$t_2\in H^{0,0}$.  This uses remarks (e) and (f) from
Section~\ref{se:back}, together with the fact that $(\Sq^1 s_1)s_2 =
\Sq^1(s_1s_2)$ for $s_2\in H^{2,2}$ (by the Cartan formula).

\vspace{0.2in} 

The final step is to analyze the groups
$\Ext^{n-1}(H^{**},\Sigma^{n,0}H^{**})$ for $n\geq 4$; these complete
the $E^{1,*}$ column of the Adams spectral sequence.  One doesn't have
to compute them explicitly, just enough to know that every element is
decomposable as a sum of products from
$\Ext^{n-2}(H^{**},\Sigma^{n-1,0}H^{**})$ and
$\Ext^{1}(H^{**},\Sigma^{1,0}H^{**})$.

The calculations involve nothing more than what we've done so far,
except for more sweat.  It's fairly easy to write down all the cocycles
made up from the classes of types (i)-(iv) listed previously.  
All bar elements which have a $\Sq^4$ in them are cocycles, for
instance.  But note that such a bar element will either begin or end
with a $\Sq^1$ or a $\Sq^2$, so that it decomposes as a
product of smaller degree cocycles (this again depends on
\ref{se:back}(e,f)).  One also finds cocycles of the
form
\[ s[\Sq^1 | \Sq^1 | \cdots | \Sq^3 | \Sq^1 | \cdots | \Sq^1]^* +  
s[\Sq^1 | \Sq^1 | \cdots | \Sq^2\Sq^1 | \Sq^1 | \cdots | \Sq^1]^*, 
\]  
but for each of these a common $[\Sq^1]^*$ can be pulled off of either
the left or right side---again showing it to be decomposable. 

Certainly there are cocycles which are {\it not\/} decomposable, like
ones of the form
\[ s[\Sq^2 | \Sq^1 | \cdots | \Sq^1 | \Sq^3 ]^* +  
s[\Sq^2\Sq^1 | \Sq^1 | \cdots |\Sq^1 | \Sq^2]^*.
\]  
But this is the coboundary of $s[\Sq^2\Sq^1 | \Sq^1 | \cdots |
\Sq^1|\Sq^3]$, and so vanishes in $\Ext$.

Anyway, I am definitely  not going to give all the details.  But with
enough diligence one can see that all elements of
$\Ext^{n-1}(H^{**},\Sigma^{n,0}H^{**})$ for $n\geq 3$ do indeed
decompose into products.

\begin{remark}
A final note about Adem relations, for those who want to try their
hand at further calculations.  Every formula I've seen for the motivic
Adem relations---in publications or preprints---seems to either
contain typos or else is just  plain wrong.  A good test for a given formula
is to see whether it gives $\Sq^3\Sq^2=\rho\Sq^3\Sq^1$ (this formula
follows from the smaller Adem relation $\Sq^2\Sq^2=\tau \Sq^3\Sq^1$,
the derivation property of the Bockstein, the fact that $\beta^2=0$,
and the identity $\Sq^3=\beta\Sq^2$).
\end{remark}

\bibliographystyle{amsalpha}

\begin{thebibliography}{DHK}


\bibitem[AEJ]{AEJ} J. K. Arason, R. Elman, and B. Jacob, \emph{The
graded Witt ring and Galois cohomology I}, in
\emph{Quadratic and Hermetian forms}, Canadian Math. Soc. Conference
Proceedings Vol. 4 (1984),  17--50.

\bibitem[BT]{BT} H. Bass and J. Tate, \emph{The Milnor ring of a
  global field}, Algebraic $K$-theory, II: "Classical" algebraic
  $K$-theory and connections with arithmetic (Proc. Conf., Seattle,
  Wash., Battelle Memorial Inst., 1972), pp.  349--446. Lecture Notes
  in Math., Vol. 342, Springer, Berlin, 1973.


\bibitem[De]{De} A. Delzant, \emph{D\'efinition des classes de
Stiefel-Whitney d'un module quadratique sur un corps de
charact\'eristique diff\'erent de 2}, C.R. Acad. Sci. Paris {\bf 255},
1366--1368.

\bibitem[Du]{Du} D. Dugger, \emph{An Atiyah-Hirzebruch spectral
sequence for $KR$-theory}, to appear in $K$-theory.

\bibitem[EL]{EL} R. Elman and T.Y. Lam, \emph{Pfister forms and
$K$-theory of fields}, Jour. of Algebra {\bf 23} (1972), 181--213.


\bibitem[Ka]{Ka} K. Kato, \emph{A generalization of local class field theory 
by using $K$-groups II}, J. Fac. Sci. Univ. Tokyo Sect. 1A Math. {\bf
  27} (1980), no. 3, 603--683.

\bibitem[LZ]{LZ} J. Lannes and S. Zarati, \emph{Invariants de Hopf
d'ordre sup\'erieur et suite spectrale d'Adams}, C.R. Acad. Sc. Paris
{\bf 296} (1983), 695--698.

\bibitem[MVW]{MVW} C. Mazza, V. Voevodsky, and C. Weibel, \emph{Lectures on
motivic cohomology}, preprint, July
2002. http://www.math.uiuc.edu/K-theory/0486. 

\bibitem[M]{M} A. Merkujev, \emph{On the norm residue symbol of degree
$2$}, Soviet Math. Doklady {\bf 24} (1981), 546--551.

\bibitem[Mi]{Mi} J. Milnor, \emph{Algebraic $K$-theory and quadratic
forms}, Invent. Math. {\bf 9} (1970), 318--344.

\bibitem[M1]{M1} F. Morel, \emph{Voevodsky's proof of Milnor's
conjecture}, Bull. Amer. Math. Soc. {\bf 35}, no. 2 (1998), 123--143.

\bibitem[M2]{M2} F. Morel, \emph{Suite spectral d'Adams et invariants
cohomologiques des formes quadratiques}, C.R. Acad. Sci. Ser. 1
Math. {\bf 328} (1999), no. 11, 963--968.


\bibitem[M3]{M3} F. Morel, \emph{An introduction to $\A^1$-homotopy
theory}, Trieste lectures.  Preprint, 2002.  Available at
http://www.math.jussieu.fr/$\sim$morel/.

\bibitem[M4]{M4} F. Morel, \emph{Milnor's conjecture on quadratic
forms and mod $2$ motivic complexes}, preprint, 2004.
http://www.math.uiuc.edu/K-theory/0684. 

\bibitem[M5]{M5} F. Morel, \emph{Suites spectrales d'Adams et
conjectures de Milnor}.  Draft available online at
http://www.math.jussieu.fr/$\sim$morel/.

\bibitem[OVV]{OVV} D. Orlov, A. Vishik, and V. Voevodsky, \emph{An
exact sequence for $K^M_*/2$ with applications to quadratic forms},
preprint, 2000.  http://www.math.uiuc.edu/K-theory/0454.

\bibitem[Pf1]{Pf1} A. Pfister, \emph{Some remarks on the historical
development of the algebraic theory of quadratic forms}, in
\emph{Quadratic and Hermetian forms}, Canadian Math. Soc. Conference
Proceedings Vol. 4 (1984),  1--16.

\bibitem[Pf2]{Pf2} A. Pfister, \emph{On the Milnor conjectures:
history, influence, applications}, Jarhes. Deutsch. Math.-Verein. {\bf
102} (2000), 15--41.

\bibitem[R1]{R1} M. Rost, \emph{Some new results on the Chow groups of
  quadrics}, preprint, 1990.  Available at 
http://www.math.uiuc.edu/K-theory/0165/.

\bibitem[R2]{R2} M. Rost, \emph{Norm varieties and algebraic
cobordism}, Proceedings of the International Congress of
Mathematicians, Vol. II (Beijing, 2002), 77--85, Higher Ed. Press,
Beijing, 2002.

\bibitem[S1]{S1} W. Scharlau, \emph{Quadratic and Hermitian forms},
Grundlehren der mathematischen Wissenschaften {\bf 270},
Springer-Verlag Berlin Heidelberg, 1985.

\bibitem[S2]{S2} W. Scharlau, \emph{ On the history of the algebraic
theory of quadratic forms}, in \emph{Quadratic forms and their
applications}, Contemp. Math. {\bf 272}, American Mathematical
Society, 2000,   229--259.

\bibitem[Se]{Se} J.-P. Serre, \emph{Cohomologie Galoisienne, Cinqui\`eme
\'edition}, Lecture Notes in Math. {\bf 5}, Springer-Verlag Berlin
Heidelberg, 1973, 1994.

\bibitem[Su]{Su} A. Suslin, \emph{Voevodsky's proof of the Milnor
conjecture}, Current Developments in Mathematics, 1997 (Cambridge,
MA), 173--188.

\bibitem[V1]{Vmilnor} V. Voevodsky, \emph{The Milnor conjecture}, preprint,
1996.   http://www.math.uiuc.edu/K-theory/0170.

\bibitem[V2]{V2} V. Voevodsky, \emph{Reduced power operations in motivic
cohomology}, Publ. Math. Inst. Hautes \'Etudes Sci., No. 98
(2003), 1-57.

\bibitem[V3]{Vmilnor2} V. Voevodsky, \emph{Motivic cohomology with
$\Z/2$-coefficients}, Publ. Math. Inst. Hautes \'Etudes Sci., No. 98
(2003), 59--104.

\bibitem[V5]{Vbkato} V. Voevodsky, \emph{On motivic cohomology with
$\Z/l$ coefficients}, preprint, 2003.  Available at
http://www.math.uiuc.eud/K-theory/0639.

\bibitem[VSF]{VSF} V. Voevodsky, A. Suslin, and E. M. Friedlander, 
\emph{Cycles, transfers, and motivic homology theories}, 
Annals of Mathematics Studies {\bf 143}, Princeton University Press, 
Princeton, NJ, 2000.

\end{thebibliography}

\end{document}